\documentclass[11pt]{article}

\usepackage[mathscr]{eucal}%
\usepackage{amssymb}%
\usepackage{amsmath}%
\usepackage{amscd}%

\usepackage{indentfirst}%
\usepackage{amsthm}%

\newtheoremstyle{cute}
    {} 
    {} 
    {\slshape} 
    {} 
    {\bfseries} 
    { } 
    { } 
    {} 

\newtheoremstyle{definitioncute}
    {} 
    {} 
    {\upshape} 
    {} 
    {\bfseries} 
    { } 
    { } 
    {} 

\theoremstyle{cute}
\newtheorem{theorem}{Theorem}[section]
\newtheorem{lemma}[theorem]{Lemma}
\newtheorem{proposition}[theorem]{Proposition}
\newtheorem{corollary}[theorem]{Corollary}

\newtheorem*{namedtheorem}{\theoremname}
\newcommand{\theoremname}{Theorem}
\newenvironment{named}[1]{\renewcommand{\theoremname}{#1}\begin{namedtheorem}}%
    {\end{namedtheorem}}

\theoremstyle{definitioncute}
\theoremstyle{remark}
\newtheorem*{remark}{Remark}
\newtheorem*{remarks}{Remarks}

\newcounter{chislo}
\newenvironment{myenumerate}%
    {\begin{list}{$\mathrm{(\roman{chislo})}$}{\usecounter{chislo}%
        \setlength{\topsep}{.5ex}\setlength{\labelwidth}{50pt}%
        \setlength{\parsep}{0pt}\setlength{\itemsep}{0.5ex}%
        \setlength{\leftmargin}{27pt}}}%
    {\end{list}}
\makeatletter

\renewcommand{\@seccntformat}[1]{\csname the#1\endcsname.\hspace{.5em}}

\renewcommand{\section}{\@startsection
    {section}
    {1}
    {0em}
    {1.5\baselineskip plus .5\baselineskip minus .5\baselineskip}
    {0.5\baselineskip plus .2\baselineskip minus 0\baselineskip}
    {\normalfont\normalsize\scshape\large\centering}}

\renewcommand{\subsection}{\@startsection
    {subsection}
    {2}
    {0em}
    {.5\baselineskip plus .5\baselineskip minus 0\baselineskip}
    {.5\baselineskip}
    {\normalfont\normalsize\scshape}}

\makeatother

\numberwithin{equation}{section}

\renewcommand{\title}[1]{\begin{center}\bfseries\Large #1\end{center}\smallskip}
\renewcommand{\author}[1]{\begin{center}\scshape\normalsize #1\end{center}}
\newcommand{\oo}{\ensuremath{\mathcal{O}}}
\newcommand{\cali}[1]{\ensuremath{\mathcal{#1}_v}}
\newcommand{\congruence}[3]{\ensuremath{#1 \equiv #2 \ ( \mathrm{mod} \ #3)}}
\newcommand{\tg}{\ensuremath{G}}
\DeclareMathOperator{\Ker}{\ensuremath{\mathrm{Ker}}}%
\DeclareMathOperator{\Image}{\ensuremath{\mathrm{Im}}}
\begin{document}

\title{Bounded Generation of $S$-arithmetic Subgroups of Isotropic Orthogonal
Groups over Number Fields}

\author{Igor V.\ Erovenko and Andrei S.\ Rapinchuk}

\begin{abstract}
Let $f$ be a nondegenerate quadratic form in $n \geqslant 5$ variables over a
number field $K$ and let $S$ be a finite set of valuations of $K$ containing
all Archimedean ones. We prove that if the Witt index of $f$ is $\geqslant 2$
or it is $1$ and $S$ contains a non-Archimedean valuation, then the
$S$-arithmetic subgroups of $\mathbf{SO}_n(f)$ have bounded generation. These
groups provide a series of examples of boundedly generated $S$-arithmetic
groups in isotropic, but not quasi-split, algebraic groups.
\end{abstract}

\smallskip

\centerline{\large\it In memory of Professor Gordon E.\ Keller}

\smallskip

\section{Introduction}

An abstract group $\Gamma$ is said to have \emph{bounded generation} 
(abbreviated (BG)) if there exist elements $\gamma_1, \ldots , \gamma_t \in 
\Gamma$ such that $\Gamma = \langle \gamma_1 \rangle \cdots \langle \gamma_t 
\rangle$, where $\langle \gamma_i \rangle$ is the cyclic subgroup generated by 
$\gamma_i$. Such groups are known to have a number of remarkable properties: 
the pro-$p$ completion $\hat{\Gamma}^{(p)}$ is a $p$-adic analytic group for 
every prime $p$ \cite{DdSMS,lazard};  if $\Gamma$ in addition satisfies 
condition (Fab)\footnote{We recall that condition (Fab) for $\Gamma$ means 
that every subgroup of finite index $\Gamma_1$ of $\Gamma$ has finite 
abelianization $\Gamma_1^{\rm ab} = \Gamma_1 / [\Gamma_1 , \Gamma_1].$} then 
it has only finitely many inequivalent completely reducible representations in 
every dimension $n$ over any field (see \cite{lubotzky-86,R2,rapinchuk-98} for 
representations in characteristic zero, and \cite{ALuP} for arbitrary 
characteristic); if $\Gamma$ is an $S$-arithmetic subgroup of an absolutely 
simple simply connected algebraic group over a number field, then $\Gamma$ has 
the congruence subgroup property \cite{lubotzky-95,platonov-rapinchuk-92}. 
There are reasons to believe that the class of groups having (BG) is 
sufficiently broad, in particular it most probably contains all higher rank 
lattices in characteristic zero (note that there are also \emph{simple,} hence 
nonlinear, infinite boundedly generated groups \cite{Mur}). Unfortunately, 
bounded generation of lattices is known only in very few cases. First, it was 
noted that the results on factoring a unimodular matrix over an arithmetic 
ring as a product of a bounded number of elementary matrices 
\cite{carter-keller-83,cooke-w-75,liehl-84,murty-95,carter-keller-preprint} 
imply bounded generation of the corresponding unimodular groups (notice, 
however, that the results on ``bounded factorization'' do not extend to 
``nonarithmetic'' Dedekind rings \cite{van-der-kallen}). Later, Tavgen 
\cite{tavgen-90} showed that every $S$-arithmetic subgroup of a split or 
quasi-split algebraic group over a number field $K$ of $K$-rank $\geqslant 2$ 
is boundedly generated. However, until recently there were no examples of 
boundedly generated $S$-arithmetic groups in algebraic groups that are not 
split or quasi-split. The goal of this paper is to establish bounded 
generation of a large family of $S$-arithmetic subgroups in isotropic 
orthogonal groups.

\begin{named}{Main Theorem}\label{thm:main}
Let $f$ be a nondegenerate quadratic form over a number field $K$ in $n 
\geqslant 5$ variables, $S$ be a finite set of valuations of $K$ containing all 
Archimedean ones. Assume that either the Witt index of $f$ is $\geqslant 2$ or 
it is one and $S$ contains a non-Archimedean valuation. Then any $S$-arithmetic 
subgroup of\/ $\mathbf{SO}_n(f)$ has bounded generation.
\end{named}

This result was announced  with a sketch of proof in 
\cite{erovenko-rapinchuk-01} for the case where the Witt index of $f$ is 
$\geqslant 2$. The argument in \cite{erovenko-rapinchuk-01} boiled down to 
reducing the general case to $n = 5$ where the group is split, so one can use 
the result of Tavgen \cite{tavgen-90}. Unfortunately, this argument does not 
immediately extend to the situation where the Witt index is one due to some 
technical problems, but mainly because of the fact that the resulting special 
orthogonal group in dimension $n = 5$ is no longer split and bounded 
generation of its $S$-arithmetic subgroups has not been previously 
established. At the same time, the method used in \cite{erovenko-rapinchuk-01} 
does not allow one to reduce $n = 5$ to $n = 4$ where the orthogonal group has 
type $A_1 \times A_1$ so one can apply the known results for $\mathrm{SL}_2$. 
In the present paper, the method of \cite{erovenko-rapinchuk-01} is modified 
in order to to overcome the difficulties noted above. Our primary objective 
was to treat the case $n = 5$, but it turned out that the resulting argument 
applies in all dimensions and in fact simplifies the proof given in 
\cite{erovenko-rapinchuk-01}.

Now, we explain briefly how the proof of the Main Theorem goes. To facilitate 
the use of strong approximation, we will argue for the spin group $G = {\bf 
Spin}_n(f)$ rather than for the special orthogonal group ${\bf SO}_n(f);$ 
notice that (BG) of $S$-arithmetic subgroups in one of them implies the same 
property for the other --- see Proposition~\ref{prop:isogeny}. We consider the 
standard representation of $G$ on the $n$-dimensional quadratic space and, 
after choosing appropriately two anisotropic orthogonal vectors $a$, $b \in 
K^n$, analyze the product map
\[
P := G(a) \times G(b) \times G(a) \times G(b) \stackrel{\mu}{\longrightarrow} 
G,
\]
where $G(a)$ and $G(b)$ denote the stabilizers of $a$ and $b$, respectively. 
The proof of (BG) is reduced from dimension $n$ to dimension $n - 1$ by 
proving that either $\mu(P_{\oo(S)})$ or a product of its several copies 
contains a subset of $G_{\oo(S)}$ open with respect to the topology defined by 
a certain finite set of valuations (see \S\ref{sec:proof} for precise 
formulations). To achieve this, we construct an auxiliary variety $Z$ and 
factor $\mu$ as a product of two regular maps $\phi \colon P \rightarrow Z$ 
and $\psi \colon Z \rightarrow G$, see \S\ref{sec:setup}. We then establish a 
local-global principle for the fibers of $\phi$ (see \S\ref{sec:fibers}), and 
finally make sure that the relevant local conditions are satisfied. 
Eventually, this process enables us to descend either to a 5-dimensional form 
of Witt index two or to a 4-dimensional isotropic form. So, to complete the 
argument it remains to observe that (BG) of  $S$-arithmetic subgroups is a 
result of Tavgen \cite{tavgen-90} in the first case, and follows from the 
known results for $\mathrm{SL}_2$ 
\cite{cooke-w-75,liehl-84,murty-95,carter-keller-preprint} in the second case. 
It appears that some parts of the argument, particularly the entire method of 
factoring a sizable set of $S$-integral transformations of a quadratic lattice 
as a product of transformations of sublattices having smaller rank, are of 
independent interest and may have other applications.

The first-named author would like to acknowledge partial financial support from 
the UNCG New Faculty Grant. The second-named author would like to acknowledge 
partial financial support from the NSF grant DMS-0138315, the BSF grant 
\#2000171 and the Humboldt Foundation (Bonn, Germany).

\section{Preliminaries}\label{sec:setup}

In this section, $K$ will denote an arbitrary field of characteristic $\neq 2$. 
Let $f$ be a nondegenerate quadratic form over $K$ of dimension $n \geqslant 
5$. Given an extension $E/K$, we let $i_E(f)$ denote the Witt index of $f$ over 
$E$, and we will write $i_v(f)$ instead of $i_{K_v}(f)$ if $K$ is a number 
field and $v$ is a valuation of $K$. Throughout the paper, we will assume that 
$i_K(f) \geqslant 1$, i.e., $f$ is $K$-isotropic. We realize $f$ on an 
$n$-dimensional vector space $W$ over $K$ and let $(\cdot \mid \cdot)$ denote 
the associated bilinear form. We also fix a basis $e_1, e_2, \ldots , e_n$ of 
$W$ in which $f$ looks as follows:
\begin{equation}\label{eq:basisforf}
f(x_1, \ldots , x_n) = x_1 x_2 + \alpha_3 x_3^2 + \cdots + \alpha_n x_n^2,
\end{equation}
and set $a = e_n$, $b = e_{n-1}$.

Next, we need to introduce some algebraic varieties and morphisms between them. 
Consider $\mathbf W = W \otimes_K \Omega$, where $\Omega$ is a ``universal 
domain'', and extend $f$ and $(\cdot \mid \cdot)$ to $\mathbf W$. Let $G$ 
denote the spin group $\mathbf{Spin}_n(f)$ associated with $\mathbf W$, 
regarded as an algebraic $K$-group (naturally) acting on $\mathbf W$. For a 
vector $w \in \mathbf W$, $G(w)$ will denote its stabilizer, and we will write 
$G(w_1 , w_2)$ for $G(w_1) \cap G(w_2)$ etc. We will be working with the 
following algebraic $K$-varieties:
\begin{align*}
P & = G(a) \times G(b) \times G(a) \times G(b),\\
X & = \{ s \in \mathbf{W} \mid f(s) = f(a) \},\\
Y & = \{ (g, s) \in G \times X \mid  (s | g(b)) = 0 \},\\
Z & = \{ (g,s,t) \in Y \times \mathbf W \mid (t|a) = 0,\ f(t) = f(b), \ (s|t) = 
0 \},
\end{align*}
and the following morphisms:
\begin{align*}
\mu \colon & P \to G, \quad \mu (x,y,z,u) = xyzu,\\
\phi \colon & P \to Z, \quad \phi (x,y,z,u) = (xyzu,\ xy(a),\ xy(b)),\\
\varepsilon \colon & Z \to Y, \quad \varepsilon (g, s, t) = (g, s),\\
\nu \colon & Z \to X, \quad \nu (g, s, t) = s.
\end{align*}
We notice that the image of $\phi$ is indeed contained in $Z$ as
\[
(xy(a)|xyzu(b)) = (a|zu(b)) = (z^{-1}(a)|b) = (a|b) = 0,
\]
hence $(xyzu, xy(a)) \in Y$, and also
\[
(xy(b)|a) = (b|x^{-1}(a)) = (b|a) = 0.
\]
The proof of the Main Theorem hinges on the fact that the product morphism $\mu 
\colon P \to G$ can be factored as $\mu = \psi \circ \phi$, where
\[
\psi \colon Z \to G, \quad \psi (g, s, t) = g.
\]

\begin{proposition}\label{prop:phipsi}
\
\begin{myenumerate}
\item For every $g \in G_K$, $\psi^{-1}(g)_K \neq \emptyset$.

\item For every $\zeta \in Z_K$, $\phi^{-1}(\zeta)_K \neq \emptyset$.
\end{myenumerate}
Consequently, $\mu (P_K) = G_K$.
\end{proposition}
\begin{proof}
(i) If $g(b) = \pm b,$ one easily verifies that $(g, a, b) \in \psi^{-1}(g)_K.$ 
So, we may assume that $g(b) \neq \pm b.$ Set
\[
u' = \frac{(g(b) | b)}{f(b)}\, b.
\]
Being isotropic, the space $\langle a , b \rangle^{\perp}$ contains a nonzero 
vector $u''$ such that $f(u'') = f(b) - f(u').$ Then the vector $u := u' + u''$ 
satisfies the following conditions:
\[
(u | a) = 0, \quad (u | b) = (g(b) | b), \quad \text{and} \quad f(u) = f(b).
\]
Since $g(b) \neq \pm b,$ the last two conditions imply that $\langle u, b 
\rangle$ and $\langle g(b), b \rangle$ are isometric 2-dimensional subspaces of 
$W,$ so by Witt's theorem there exists $\sigma \in \mathrm O_n(f)$ such that 
$\sigma(u) = g(b)$ and $\sigma(b) = b.$ Then
\[
(\sigma(a) | b) = (a | b) = 0
\]
and
\[
(\sigma(a) | g(b)) = (a | \sigma^{-1}(g(b))) = (a | u) = 0,
\]
implying that $(g, \sigma(a), b) \in \psi^{-1}(g)_K.$

(ii) Suppose that $\zeta =(g, s, t) \in Z_K$. Since $(t | a) = 0,$ the vectors 
$t$ and $a$ are linearly independent. As $f(t) = f(b),$ by Witt's theorem there 
exists $\rho \in \mathrm O_n(f)$ such that
\begin{equation}\label{E:1}
\rho (b) = t \quad \text{and} \quad \rho (a) = a.
\end{equation}
In fact, one can always find such a $\rho$ in $\mathrm{O}'_n(f),$ the kernel 
of the spinor norm $\theta$ on $\mathrm{SO}_n(f).$ Indeed, if $\det \rho = 
-1$, we can pick an anisotropic $c \in W$ orthogonal to both $a$ and $b$, and 
replace $\rho$ with $\rho \tau_c,$ where $\tau_c$ is the reflection associated 
with $c,$ which allows us to assume that $\rho \in \mathrm{SO}_n(f).$ 
Furthermore, since the space $\langle a, b \rangle^{\perp}$ is isotropic, 
there exists $\delta \in \mathrm{SO}_n(f)(a, b)$ such that $\theta(\delta) = 
\theta(\rho)$ (see \cite[Thm.~5.18]{artin}).  Then we can replace $\rho$ with 
$\rho\delta^{-1} \in \mathrm{O}'_n(f).$

Arguing similarly, we find $\eta \in \mathrm{O}'_n(f)$ such that
\begin{equation}\label{E:2}
\eta (a) = s \quad \text{and} \quad \eta (b) = t
\end{equation}
and $\sigma \in \mathrm{O}'_n(f)$ such that
\begin{equation}\label{E:3}
\sigma (a) = s \quad \text{and} \quad \sigma (b) = g(b).
\end{equation}
Since the elements $\rho,$ $\eta$, and $\sigma$ have spinor norm one, they are 
images under the canonical isogeny $\pi \colon \mathbf{Spin}_n(f) \rightarrow 
\mathbf{SO}_n(f)$ of suitable elements $\tilde{\rho}$, $\tilde{\eta}$, 
$\tilde{\sigma} \in G_K = \mathrm{Spin}_n(f).$ Set $x = \tilde{\rho}$, $y = 
\tilde{\rho}^{-1} \tilde{\eta}$, $z = \tilde{\eta}^{-1} \tilde{\sigma}$, and $u 
= (xyz)^{-1} g$. Then $(x, y, z, u) \in P_K$ and $xyzu = g.$ Moreover, $xy (a) 
= s$ and $xy (b) = t,$ which shows that $\phi(x, y, z, u) = \zeta,$ as 
required.
\end{proof}

\begin{remark}
It follows from Proposition~\ref{prop:phipsi} that $G_K = G(a)_K G(b)_K G(a)_K 
G(b)_K$. For classical groups, decompositions of this kind were introduced by 
M.~Borovoi \cite{borovoi-85}. Our proof of the Main Theorem is based on the 
analysis of the Borovoi decomposition for the group of $S$-integral points. In 
\cite{erovenko-rapinchuk-01} we used the Borovoi decomposition involving three 
factors, $G(a)$, $G(b)$, and $G(a)$, but as we will see, the decomposition with 
four factors allows one to bypass some technical difficulties and eventually 
leads to a more general result.
\end{remark}

The following properties of the morphisms introduced above will be used in the 
sequel.

\begin{lemma}\label{L:1A}
The morphisms $\phi \colon P \to Z$ and $\varepsilon \colon Z \to Y$ are 
surjective. Consequently, if $\mathrm{char}\, K = 0,$ there exists a Zariski 
$K$-open set $P_0 \subset P$ such that for $h \in P_0$, the points $\phi (h)$ 
and $(\varepsilon \circ \phi)(h)$ are simple on $Z$ and $Y$ respectively, and 
the differentials $d_h \phi \colon T(P)_h \to T(Z)_{\phi (h)}$, $d_{\phi (h)} 
\varepsilon \colon T(Z)_{\phi (h)} \to T(Y)_{(\varepsilon \circ \phi) (h)}$, 
and $d_h \mu \colon T(P)_h \to T(G)_{\mu (h)}$ are surjective.
\end{lemma}
\begin{proof}
It follows from Proposition~\ref{prop:phipsi} that $\phi \colon P \to Z$ and 
$\mu \colon P \to G$ are surjective. Now, given $(g ,s) \in Y$, over an 
algebraically closed field one  can always find $t \in \langle a, s 
\rangle^{\perp}$ such that $f(t) = f(b)$, whence the  surjectivity of 
$\varepsilon \colon Z \to Y$. The rest of the lemma follows from a well-known 
result about dominant separable morphisms (see, for example, \cite[Ch.~AG, 
Thm.~17.3]{Borel}) and the irreducibility of $P.$
\end{proof}

\begin{lemma}\label{L:2A}
Set $\eta = \nu \circ \phi$. Then $\eta(P_E) = X_E$ for any extension $E/K$.
\end{lemma}
\begin{proof}
It is enough to show that $\phi (P_E) = Z_E$ and $\nu (Z_E) = X_E,$ the first 
assertion being part (ii) of Proposition~\ref{prop:phipsi} in which $K$ is 
replaced with $E.$ For the second assertion, let $s \in X_E.$ Then by Witt's 
theorem there exists $g \in \mathbf{SO}_n(f)_E$ such that $g(a) = s.$ Since the 
orthogonal complement to $a$ in $W \otimes_K E$ is isotropic, arguing as in the 
proof of part (ii) of Proposition~\ref{prop:phipsi}, we see that $g$ can be 
chosen to be of the form $g = \pi(\tilde{g})$ for some $\tilde{g} \in G_E.$ If 
$s = \pm a$, then one immediately verifies that $(\tilde{g}, s, b) \in 
\nu^{-1}(s)_E.$ Otherwise, the space $\langle a, s \rangle$ is 2-dimensional. 
Since the orthogonal complement to $\langle a, b \rangle$ in $W \otimes_K E$ 
is isotropic, we can argue as in the proof of part (i) of 
Proposition~\ref{prop:phipsi} to find $w \in W \otimes_K E,$ $w \notin \langle 
a , b \rangle$, satisfying
\[
(w | b) = 0, \quad (w | a ) = (s | a), \quad \text{and} \quad f(w) = f(a).
\]
By Witt's theorem, it follows from the last two conditions that there exists 
$\sigma \in \mathbf{SO}_n(f)_E$ such that $\sigma(w) = s$ and $\sigma(a) = a.$ 
Set $t = \sigma(b) \in W \otimes_K E.$ Then
\[
(t | a) = (\sigma(b) | \sigma(a)) = (b | a) = 0
\]
and
\[
(t | s) = (\sigma (b) | g(a)) = (b | \sigma^{-1} g (a)) = (b | w) = 0
\]
whence $(g, s, t) \in \nu^{-1}(s)_K$, as required.
\end{proof}

\section{Fibers of the morphism $\phi$}\label{sec:fibers}

From now on,  $K$ will denote a number field. We let  $V^K$, $V^K_{\infty}$, 
and $V^K_f$ denote the set of all inequivalent valuations of $K$, the subsets 
of Archimedean, and non-Archimedean valuations, respectively. As usual, for $v 
\in V^K$, $K_v$ denotes the completion of $K$ with respect to $v$, and for $v 
\in V^K_f,$ $\mathcal O_v$ denotes the ring of integers in $K_v$ (by 
convention, $\mathcal O_v = K_v$ for $v \in V^K_{\infty}$). Given a finite 
subset $S$ of $V^K$ containing $V^K_{\infty}$, we let ${\mathcal O}(S)$ denote 
the ring of $S$-integers in $K,$ i.e.,
\[
\mathcal O(S) = \{ x \in K \mid x \in \mathcal O_v \ \text{for all} \ v \notin 
S \}.
\]
Finally, $A_{K, S}$ will denote the ring of $S$-adeles of $K$ (adeles without 
the components corresponding to the valuations in $S$), and $A_{K, S}(S) := 
\prod_{v \notin S} {\mathcal O}_v$ will be the ring of $S$-integral $S$-adeles.

Now, suppose that $f$ is a quadratic form as in \S\ref{sec:setup}. For a real 
$v \in V^K_{\infty}$, we let $(n_v^+ , n_v^-)$ denote the signature of $f$ over 
$K_v = \mathbb R$. By scaling $f$ (which does not affect the orthogonal group), 
we can achieve that $n_v^+ \geqslant n_v^-$ (and consequently $n_v^+ \geqslant 
3$ as $n \geqslant 5$) for all real $v \in V^K_{\infty}$. Then one can choose a 
basis $e_1, \ldots , e_n$ of $W = K^n$ so that in the corresponding expression 
\eqref{eq:basisforf} for $f$, the coefficients  $\alpha_i$ belong to $\mathcal 
O(S)$ for all $i = 3, \ldots , n$, and, in addition, $\alpha_{n-1}$, $\alpha_n
> 0$ in $K_v = \mathbb R$ for all real $v \in V^K_{\infty}$ (these conventions
will be kept throughout the rest of the paper).

As we mentioned in the previous section, our goal is to find a version of the 
Borovoi decomposition for the group of $S$-integral points. Towards this end, 
in this section we will develop some conditions on $\zeta \in Z_{\mathcal 
O(S)}$ which ensure that $\phi^{-1}(\zeta)_{\mathcal O(S)} \neq \emptyset$. To 
avoid any ambiguity, we would like to stipulate that $S$-integral points in the 
space $\mathbf{W}$ and its (closed) subvarieties will be understood relative to 
the fixed basis $e_1, \ldots , e_n,$ and $G_{{\mathcal O}(S)}$ by definition 
consists of those $g \in G_K$ for which $\pi(g) \in \mathrm{SO}_n(f)$ is 
represented in the fixed basis by a matrix with entries in $\oo(S)$ (of course, 
it is possible to realize $G_{\oo(S)}$ as the group of $S$-integral points in 
the usual sense with respect to some faithful representation of $G,$ but we 
will not need this realization). The same conventions apply to ${\mathcal 
O}_v$-points for $v \notin S.$

The following set of valuations plays a prominent role in our argument:
\[
V_0 = \left( \cup_{i=3}^n V(\alpha_i) \right) \cup V(2),
\]
where for $\alpha \in K^{\times}$, we set $V(\alpha) = \{ v \in V^K \setminus S 
\mid v(\alpha) \neq 0 \}$.

\begin{theorem}\label{prop:application-witt}
Let $\zeta \in Z_{\oo(S)}$. Suppose that $\phi^{-1}(\zeta)_{\oo_v} \neq 
\emptyset$ for all $v \in V_0$. Then $\phi^{-1}(\zeta)_{\oo(S)} \neq 
\emptyset$.
\end{theorem}

We begin by establishing the following local-global principle for the fibers of 
$\phi$.

\begin{lemma}\label{prop:localglobal}
Let $\zeta \in {Z}_{\oo(S)}$. Suppose that ${\phi}^{-1}(\zeta)_K \neq 
\emptyset$ and ${\phi}^{-1}(\zeta)_{\oo_v} \neq \emptyset$ for all $v \notin 
S$. Then ${\phi}^{-1}(\zeta)_{\oo(S)} \neq \emptyset$.
\end{lemma}
\begin{proof}
Let $\zeta = (g,s,t)$ and $H  =  G(a,b)$. Being the spinor group of the space 
$\langle a , b \rangle^{\perp},$ which  is $K$-isotropic and has dimension 
$\geqslant 3,$ the group $H$ has the property of strong approximation with 
respect to $S,$ i.e., (diagonally embedded) $H_K$ is dense in $H_{A_{K , S}}$ 
(see \cite[104:4]{omeara}, \cite[Thm.~7.12]{platonov-rapinchuk-book}). We now 
observe that $\phi^{-1}(\zeta)$ is a principal homogeneous space of the group 
$H \times H \times H.$ More precisely, the equation
\begin{equation}\label{E:action}
(h_1, h_2, h_3) \cdot (x, y, z, u) = (x h_1^{-1}, h_1 y h_2^{-1}, h_2 z 
h_3^{-1}, h_3u)
\end{equation}
defines a simply transitive action of $H \times H \times H$ on 
$\phi^{-1}(\zeta).$ Indeed, one immediately verifies that for any $(x , y, z , 
u) \in \phi^{-1}(\zeta)$ and any $(h_1 , h_2 , h_3) \in H \times H \times H,$ 
the right-hand side of (\ref{E:action}) belongs to $\phi^{-1}(\zeta),$ and that 
$(\ref{E:action})$ defines an action. Now, suppose that $(x_i , y_i , z_i , 
u_i) \in \phi^{-1}(\zeta),$ where $i = 1 , 2.$ Set
\[
h_1 = x_2^{-1} x_1, \quad h_2 = (x_2 y_2)^{-1} (x_1 y_1), \quad h_3 = (x_2 y_2 
z_2)^{-1} (x_1 y_1 z_1).
\]
Then the conditions $x_i(a) = a$ and $x_i(b) = t$ for $i = 1, 2$ imply that 
$h_1 \in H.$ Similarly, from $(x_i y_i)(a) = s$ and $(x_i y_i)(b) = t$ we 
derive that $h_2 \in H,$ and from $(x_i y_i z_i)(a) = s$ and $(x_i y_i z_i)(b) 
= g(b)$ that $h_3 \in H.$ In view of our construction, to prove that
\[
(x, y, z, u) := (h_1, h_2, h_3) \cdot (x_1, y_1, z_1, u_1)
\]
coincides with $(x_2, y_2, z_2, u_2),$ it remains to observe that
\[
u = h_3 u_1 = (x_2 y_2 z_2)^{-1} (x_1 y_1 z_1) u_1 = (x_2 y_2 z_2)^{-1} g = 
u_2,
\]
so our claim follows.

Now, fix $(x, y, z, u) \in \phi^{-1}(\zeta)_K.$ Then
\[
\Sigma = \{ (h_1, h_2, h_3) \in  H_{A_{K, S}} \times H_{A_{K, S}} \times 
H_{A_{K, S}} \ \mid \ (h_1, h_2, h_3) \cdot (x, y, z, u) \in 
\phi^{-1}(\zeta)_{A_{K, S}(S)} \}
\]
is a nonempty open subset of $H_{A_{K , S}} \times H_{A_{K , S}} \times H_{A_{K 
, S}}.$ By strong approximation for $H,$ there exists $(h_1, h_2, h_3) \in (H_K 
\times H_K \times H_K) \cap \Sigma,$ and then
\[
(h_1, h_2, h_3) \cdot (x, y, z, u) \in \phi^{-1}(\zeta)_K \cap 
\phi^{-1}(\zeta)_{A_{K, S}(S)} = \phi^{-1}(\zeta)_{\oo(S)}
\]
is a required $S$-integral point.
\end{proof}

To finish the proof of Theorem~\ref{prop:application-witt} it now remains to 
prove the following.

\begin{lemma}\label{lm:localanalog}
Let $\zeta \in Z_{\oo(S)}$. Then $\phi^{-1}(\zeta)_{\oo_v} \neq \emptyset$ for 
all $v \notin S \cup V_0$.
\end{lemma}

The proof of Lemma~\ref{lm:localanalog} requires a version of Witt's theorem 
for local lattices, which we will state now and prove in the next section. Fix 
$v \in V^K_f,$ and let $w_1, \ldots, w_n$ be an arbitrary basis of $W_v = W 
\otimes_K K_v$ in which the matrix $F$ of the quadratic form $f$ has entries 
in ${\mathcal O}_v.$ Consider the $\mathcal{O}_v$-lattice $L_v$ with the basis 
$w_1, \ldots , w_n,$ its reduction $\bar{L}^{(v)} = L_v/\mathfrak p_v L_v$ 
modulo $\mathfrak{p}_v$ (which is an $n$-dimensional vector space over $k_v = 
{\mathcal O}_v / \mathfrak{p}_v$) and the corresponding reduction map $L_v \to 
\bar{L}^{(v)},$ $l \mapsto \bar{l}.$ We also let
\[
\mathrm O_n(f)^{L_v}_{\mathcal{O}_v} = \{ \sigma \in \mathrm O_n(f) \mid \sigma 
(L_v) = L_v \}
\]
be the stabilizer of $L_v.$

\begin{theorem}[Witt's theorem for local lattices]\label{thm:wittorthogonal}
Suppose that the systems $\{ a_1, \ldots , a_m \}$ and $\{ b_1, \ldots , b_m 
\}$ of vectors in $L_v$ satisfy the following properties:
\begin{myenumerate}
\item $(a_i | a_j)=(b_i | b_j)$ for all $1 \leqslant i \leqslant j \leqslant 
m;$

\item the systems $\{ \bar{a}_1, \ldots , \bar{a}_m \}$ and $\{ \bar{b}_1, 
\ldots , \bar{b}_m \}$ obtained  by reduction modulo $\mathfrak{p}_v$ are both 
linearly independent over $k_v.$
\end{myenumerate}
If\/ $\det F \in \cali{O}^{\times}$ and $v(2)=0$, then there exists $\sigma \in 
\mathrm O_n(f)^{L_v}_{\mathcal{O}_v}$ such that $\sigma (a_i) = b_i$ for all $i 
= 1, \dots , m$. Moreover, if $2m+1 \leqslant n$ then such a $\sigma$ can be 
found in  $\mathrm{SO}_n(f)^{L_v}_{\mathcal{O}_v}.$
\end{theorem}

\begin{proof}[Proof of Lemma~\ref{lm:localanalog}] We mimic the proof of
Proposition~\ref{prop:phipsi}(ii) except that instead of the usual Witt's 
theorem we use Theorem~\ref{thm:wittorthogonal}. We let $L$ denote the 
${\mathcal O}(S)$-lattice with the basis $e_1, \ldots , e_n$ which was fixed 
earlier, and for $v \notin S$ we set $L_v = L \otimes_{\oo(S)} \oo_v.$ We claim 
that for every $v \in V^K \setminus (S \cup V_0),$ each of the following three 
pairs $(\overline{t}^{(v)}, \overline{a}^{(v)}),$ $(\overline{s}^{(v)}, 
\overline{t}^{(v)}),$ and $(\overline{s}^{(v)}, \overline{g(b)}^{(v)})$ (where 
$\bar{\ }^{(v)}$ denotes the reduction map modulo $\mathfrak{p}_v$) is linearly 
independent over $k_v.$ For this we notice that $f(s) = f(a) = \alpha_n$ and 
$f(t) = f(b) = \alpha_{n-1}$ (cf.\ \eqref{eq:basisforf}), and because $v \notin 
V_0 \cup S,$ both $\alpha_{n-1}$ and $\alpha_n$ are invertible in $\oo_v.$ Now, 
if for example, $\overline{t}^{(v)} = \lambda \overline{a}^{(v)}$ with $\lambda 
\in k_v,$ then the condition $(t | a) = 0$ implies that
\[
\bar{0} = (\overline{t}^{(v)} | \overline{a}^{(v)}) = \lambda 
(\overline{a}^{(v)} | \overline{a}^{(v)}) = \lambda \bar{\alpha}_n.
\]
So, $\lambda = 0.$ But $\overline{t}^{(v)} \neq \bar{0}$ as $\bar{\alpha}_{n-1} 
\neq \bar{0},$ a contradiction. All other cases are considered similarly using 
the orthogonality relations in the definition of $Z.$ Then using 
Theorem~\ref{thm:wittorthogonal} we find elements $\rho,$  $\eta,$ and 
$\sigma$ in $\mathrm{SO}_n(f)^{L_v}_{\oo_v}$ satisfying conditions 
\eqref{E:1}, \eqref{E:2}, and \eqref{E:3}. Since $v \notin V_0 \cup S,$ the 
lattice $L_v$ is unimodular, and therefore the spinor norm of all three 
elements belongs to $\mathcal{O}^{\times}_v K_v^{\times 2}$ 
\cite[92:5]{omeara}. On the other hand, the lattice $M_v := L_v \cap \langle 
a, b \rangle^{\perp}$ (which has $e_1, \ldots , e_{n-2}$ as its $\oo_v$-basis) 
is unimodular of rank $\geqslant 3,$ implying that $\theta(\mathrm{SO}_n(f)(a, 
b)^{M_v}_{\oo_v}) = \mathcal{O}^{\times}_v K_v^{\times 2}$ 
\cite[92:5]{omeara}. So, arguing as in the proof of 
Proposition~\ref{prop:phipsi}(ii), we can modify the elements $\rho,$ $\eta,$ 
and $\sigma$ so that they all have trivial spinor norm. Then they can be 
lifted to elements $\tilde{\rho},$ $\tilde{\eta},$ and $\tilde{\sigma}$ in 
$G_{\oo_v},$ and one easily verifies that the quadruple $(x, y, z, u),$ where 
$x = \tilde{\rho}$, $y = \tilde{\rho}^{-1} \tilde{\eta}$, $z = 
\tilde{\eta}^{-1} \tilde{\sigma}$, and $u = (xyz)^{-1} g,$ belongs to 
$\phi^{-1}(\zeta)_{\oo_v},$ proving the lemma.
\end{proof}

The proof of Theorem~\ref{prop:application-witt} is now complete.

\section{Proof of Witt's theorem for local lattices}\label{sec:witt}

We will prove Theorem~\ref{thm:wittorthogonal} in a more general situation 
then that we dealt with in the previous section. Namely, let $\mathscr{K}$ be 
a field of characteristic $\neq 2$ which is complete with respect to a 
discrete valuation $v$ (see a remark at the end of the section regarding 
generalizations to not necessarily complete discretely valued fields). We let 
$\mathscr{O}$, $\mathfrak{p}$, and $k$ denote the corresponding valuation ring, 
the valuation ideal, and the residue field, respectively; we also pick a 
uniformizer $\pi \in \mathfrak{p}$ so that $\mathfrak{p} = \pi\mathscr{O}.$ 
Furthermore, let $\mathscr{W}$ be an $n$-dimensional vector space over 
$\mathscr{K},$ and $f$ be an \emph{arbitrary} quadratic form on $\mathscr{W}$ 
(in particular, we are not assuming that $f$ has form \eqref{eq:basisforf} in 
a suitable basis of $\mathscr{W}$) with associated symmetric bilinear form 
$(\cdot \mid \cdot).$ We fix a basis $w_1, \ldots , w_n$ of $\mathscr{W}$ such 
that $(w_i | w_j) \in \mathscr{O}$ for all $i, j = 1, \ldots , n,$ and let 
$\mathscr{L}$ denote the $\mathscr{O}$-lattice $\mathscr{O} w_1 + \cdots + 
\mathscr{O} w_n.$ Let $\bar{\mathscr{L}} = \mathscr{L}/{\mathfrak p}{\mathscr 
L}$ be the reduction of $\mathscr{L}$ modulo $\mathfrak{p},$ with the 
corresponding reduction map $\mathscr{L} \to \bar{\mathscr{L}},$ $l \mapsto 
\bar{l}.$ In the sequel, matrix representations for linear transformations of 
$\mathscr{W}$ will be considered exclusively relative to the basis $w_1, 
\ldots , w_n;$ in particular, the $\mathscr{O}$-points of the orthogonal group 
$\mathrm{O}_n(f)$ are described as follows:
\[
\mathrm O_n(f)_{\mathscr{O}} = \{ X \in \mathrm{GL}_n(\mathscr{O}) \mid 
\hskip1pt ^{t}\!X F X = F \},
\]
where $F = ((w_i | w_j))$ is the Gram--Schmidt matrix of the form $f.$ The 
following two assumptions will be kept throughout the section:
\begin{myenumerate}
\item $\mathscr{L}$ is unimodular, i.e., $\det F \in \mathscr{O}^{\times}$;

\item $\mathrm{char}\, k \neq 2$.
\end{myenumerate}
Under these assumptions, we will prove the following, which in particular 
yields Theorem~\ref{thm:wittorthogonal}:
\begin{itemize}
\item[$(*)$] given two systems of vectors $\{ a_1, \ldots , a_m \}$ and $\{ 
b_1, \ldots , b_m \}$ in $\mathscr{L}$ satisfying conditions (i) and (ii) of 
Theorem~\ref{thm:wittorthogonal}, i.e., $(a_i | a_j) = (b_i | b_j)$ for all $i 
, j = 1, \ldots , m$, and the reduced systems $\{ \bar{a}_1, \ldots , 
\bar{a}_m \}$ and $\{ \bar{b}_1, \ldots , \bar{b}_m \}$ are linearly 
independent over $k,$ there exists $X \in \mathrm O_n(f)_{\mathscr{O}}$ with 
the property $X a_i = b_i$ for all $i = 1, \ldots , m,$ and, moreover, if $2m 
+ 1 \leqslant n,$ then such an $X$ can already be found in $\mathrm{SO}_n 
(f)_{\mathscr{O}}$.
\end{itemize}
The proof uses the standard approximation procedure due to Hensel, although we 
bypass a direct usage of Hensel's Lemma for algebraic varieties. We begin with 
a couple of lemmas.

\begin{lemma}\label{lm:upbyone}
Given an integer $l \geqslant 1$ and a matrix $X \in M_n(\mathscr{O})$ 
satisfying
\begin{equation}\label{eqn:modl}
\congruence{^tXFX}{F}{\mathfrak{p}^l},
\end{equation}
there exists $Y \in M_n(\mathscr{O})$ such that
\begin{equation}\label{eqn:modlplusone}
\congruence{^tYFY}{F}{\mathfrak{p}^{l+1}}
\end{equation}
and
\[
\congruence{Y}{X}{\mathfrak{p}^l}.
\]
\end{lemma}
\begin{proof}
We need  to find $Z \in M_n(\mathscr{O})$ for which $Y := X + \pi^l Z$ 
satisfies \eqref{eqn:modlplusone}. According to \eqref{eqn:modl}, $F - \, 
^tXFX = \pi^l A$, for some (necessarily symmetric) matrix $A \in 
M_n(\mathscr{O}).$ In view of the congruence
\[
\congruence{^tYFY}{\, ^tXFX + \pi_v^l  (\, ^tZFX + \, 
^tXFZ)}{\mathfrak{p}^{l+1}},
\]
to satisfy \eqref{eqn:modlplusone} it is enough choose $Z \in 
M_n(\mathscr{O})$ so that
\[
\congruence{^tZFX + \, ^tXFZ}{A}{\mathfrak{p}}.
\]
However, it follows from our assumptions that
\[
Z := \frac{^{t}(FX)^{-1}A}{2} \in M_n(\mathscr{O})
\]
and moreover
\[
^tZFX+ \, ^tXFZ = \, ^tZ(FX) + \, ^t(FX)Z = \hskip0pt 
^{t}\!\left(\frac{A}{2}\right) + \left(\frac{A}{2}\right) = A
\]
as $A$ is symmetric. Thus, $Z$ is as required.
\end{proof}

\begin{corollary}\label{cor:hensels}
Notations as in Lemma~\ref{lm:upbyone}, there exists $\hat{X} \in \mathrm 
O_n(f)_{\mathscr{O}}$ satisfying $\congruence{\hat{X}}{X}{\mathfrak{p}^l}.$
\end{corollary}
\begin{proof}
Using Lemma~\ref{lm:upbyone}, we construct a sequence of matrices $X_i \in 
M_n(\mathscr{O})$, $i=l, l+1 , \ldots$ such that $X_l = X$, 
$\congruence{^tX_iFX_i}{F}{\mathfrak{p}^i}$, and 
$\congruence{X_{i+1}}{X_i}{\mathfrak p^i}$ for all $i \geqslant l$. Then 
$\congruence{X_i}{X_j} {\mathfrak{p}^j}$ for all $i \geqslant j \geqslant l,$ 
implying that $\{ X_i \}$ is a Cauchy sequence in $X + M_n(\mathfrak{p}^l) 
\subset M_n(\mathscr{K}).$ As $\mathscr{K}$ is complete and $\mathfrak{p}^l$ 
is closed in $\mathscr{K},$ this sequence converges to some $\hat{X} \in X + 
M_n (\mathfrak{p}^l),$ which is as required.
\end{proof}

\begin{lemma}\label{lm:main}
Given an integer $l \geqslant 1$ and two systems of vectors $\{a_1, \ldots , 
a_m \}$ and $\{ b_1, \ldots , b_m \}$ in $\mathscr{L}$ as in $(*)$ satisfying
\begin{equation*}\label{eqn:congruentlattices}
\congruence{a_i}{b_i}{\mathfrak{p}^l} \quad \text{for all} \ i = 1, \dots, m,
\end{equation*}
there exists $X \in M_n(\mathscr{O})$ such that
\begin{gather*}
\congruence{X}{E_n}{\mathfrak{p}^l},\\
\congruence{^tXFX}{F}{\mathfrak{p}^{l+1}},
\end{gather*}
and
\[
\congruence{X a_i}{b_i}{\mathfrak{p}^{l+1}} \quad \text{for all} \ i = 1, 
\dots, m.
\]
\end{lemma}
\begin{proof}
We have $b_i = a_i + \pi^l c_i$ for some $c_i \in \mathscr{L},$ and then the 
condition $(a_i | a_j) = (b_i | b_j)$ yields
\[
\congruence{(a_i|c_j)+(c_i|a_j)}{0}{\mathfrak{p}} \quad \text{for all} \ i, j 
= 1, \ldots , m.
\]
Now, suppose we can exhibit $Y \in M_n(\mathscr{O})$ such that
\begin{equation}\label{eq:y2}
\congruence{^tYF+FY}{0}{\mathfrak{p}}
\end{equation}
and
\begin{equation}\label{eq:y1}
\congruence{Y a_i}{c_i}{\mathfrak{p}} \quad \text{for all} \ i = 1, \dots, m.
\end{equation}
Then $X :=E_n+ \pi^l Y$ is as required. Indeed,
\[
^tXFX \equiv F + \pi^l (\, ^tYF+FY) \equiv F \ (\mathrm{mod} \ \mathfrak p 
^{l+1})
\]
and
\[
X a_i = a_i + \pi^l Y a_i \equiv a_i + \pi^l c_i \equiv b_i \ (\mathrm{mod} \ 
\mathfrak p^{l+1}).
\]
On the other hand, the existence of $Y$ satisfying \eqref{eq:y2} and 
\eqref{eq:y1} follows from Lemma~\ref{prop:general} below applied to the 
vector space $\mathcal{W} = \bar{\mathscr{L}}$ over the field $\mathcal{K} =  
k,$  the symmetric matrix $\mathcal{F}$ obtained by reducing $F$ modulo 
$\mathfrak{p},$ and the vectors $x_1 = \bar{a}_1, \ldots, x_m = \bar{a}_m$, 
and $y_1= \bar{c}_1, \ldots, y_m = \bar{c}_m$ in $\mathcal{W}.$
\end{proof}

\begin{lemma}\label{prop:general}
Let $\mathcal{K}$ be an arbitrary field of characteristic $\neq 2,$ and 
$\mathcal{W} = \mathcal{K}^n.$ Let $\mathcal{F}$ be a nondegenerate symmetric 
$n \times n$ matrix over $\mathcal{K},$ $(x | y) = \, ^{t}x{\mathcal F}y$ be 
the corresponding symmetric bilinear form on $\mathcal{W}$, and
\[
\mathcal{R} = \{Y \in M_n(\mathcal{K}) \mid \, ^{t}Y \mathcal F + \mathcal F Y 
= 0 \}
\]
be the corresponding space of skew-symmetric matrices. Suppose that $x_1, 
\ldots , x_m \in \mathcal{W}$ are linearly independent vectors, and set
\[
\mathcal A = \{ (y_1, \ldots , y_m) \in \mathcal{W}^m  \mid  (x_i | y_j) + 
(x_j | y_i) = 0 \ \text{for all} \ i, j = 1, \ldots , m \}
\]
and
\[
\mathcal B = \{ (Yx_1, \ldots , Yx_m)  \mid  Y \in \mathcal{R} \}.
\]
Then $\mathcal A = \mathcal B$.
\end{lemma}
\begin{proof} 
In the standard basis, matrices in $\mathcal{R}$ correspond to  linear 
operators in 
\[
\mathcal{S} = \{ Y \in \mathrm{End}_{\mathcal K}(\mathcal{W})  \mid  (Yx | y) 
+ (x | Yy) = 0 \},
\]
which in particular yields the inclusion $\mathcal B \subset \mathcal A$. So, 
it is enough to show that $\dim \mathcal A = \dim \mathcal B$.

Let $\mathcal{V}$ be the subspace of $\mathcal{W}$ spanned by $x_1, \ldots , 
x_m.$ Then $\dim \mathcal B = \dim \mathcal{S} - \dim \mathcal{T}$, where
\[
\mathcal{T} = \{ Y \in \mathcal{S} \mid Yv = 0 \ \text{for all} \ v \in 
\mathcal{V} \}.
\]
An elementary computation based on representing the transformation in 
$\mathcal{S}$ by matrices relative to a (fixed) orthogonal basis of 
$\mathcal{W}$ yields
\begin{equation}\label{eq:dimr}
\dim \mathcal{S} = \frac{n(n-1)}{2}.
\end{equation}
To calculate $\dim \mathcal{T},$ one needs to observe that $\mathcal{W}$ 
admits a basis $v_1, \ldots , v_n$ such that the vectors $v_1, \ldots , v_m$ 
form a basis of $\mathcal{V}$ and the matrix of the bilinear form $(\cdot \mid 
\cdot)$ has the following structure
\[
\left( 
\begin{array}{ccc}
& & E_r\\
& D_{n-2r} & \\
E_r & &
\end{array} 
\right)
\]
where $E_r$ is the identity matrix of size $r,$ and $D_{n-2r}$ is a 
nondegenerate diagonal matrix of size $n-2r$ (notice that $r$ is nothing but 
the dimension of the radical of $\mathcal{V},$ and in fact the vectors $v_1, 
\ldots , v_r$ form a basis of this radical). Then, by considering the matrix 
representation of transformations from $\mathcal{S}$ relative to the basis 
$v_1, \ldots , v_n,$ one finds that
\begin{equation}\label{eq:dimt}
\dim \mathcal{T} = \frac{(n-m)(n-m-1)}{2}.
\end{equation}
From \eqref{eq:dimr} and \eqref{eq:dimt} we conclude that
\begin{equation}\label{eq:ups}
\dim \mathcal B = \frac{n(n-1)}{2} - \frac{(n-m)(n-m-1)}{2}.
\end{equation}
To calculate $\dim \mathcal A$, we consider the linear functionals $f_i,$ $i = 
1, \ldots , m$, on $\mathcal{W}$ given by $f_i(w) = (w | x_i).$ Since the 
$x_i$'s are linearly independent and $\mathcal{W}$ is nondegenerate, the 
$f_i$'s are also linearly independent, implying that the linear map 
$\mathcal{W} \to \mathcal{K}^m,$ $w \mapsto (f_1(w), \ldots , f_m(w)),$ is 
surjective. It follows that the linear map $\Phi \colon \mathcal{W}^m \to 
\mathcal{K}^{m^2}$ given by
\[
\Phi (w_1, \ldots , w_m) = (f_1(w_1), \ldots , f_m(w_1), \ldots , f_1(w_m), 
\ldots , f_m(w_m))
\]
is also surjective. Let $z_{ij}$ be the coordinate in $\mathcal{K}^{m^2}$ 
corresponding to $f_i(w_j).$ Then $\mathcal A = \Phi^{-1}(\mathcal{U})$ where 
$\mathcal{U}$ is the subspace of $\mathcal{K}^{m^2}$ defined by the conditions 
$z_{ij} + z_{ji} = 0$ for all $i , j = 1, \ldots , m.$ It follows that
\[
\dim \mathcal A = mn - \mathrm{codim}_{\mathcal{K}^{m^2}}\, \mathcal{U} = mn - 
\frac{m(m+1)}{2},
\]
comparing which with \eqref{eq:ups} we obtain $\dim \mathcal A = \dim \mathcal 
B$, completing the argument.
\end{proof}

\begin{proof}[Proof of $(*)$.] 
We will inductively construct a sequence of matrices $X_s \in \mathrm 
O_n(f)_{\mathscr{O}},$ $s = 1, 2, \ldots,$  satisfying
\begin{equation}\label{E:20}
\congruence{X_t}{X_s}{\mathfrak{p}^s} \quad \text{whenever} \ t \geqslant s
\end{equation}
and
\begin{equation}\label{E:30}
\congruence{X_s a_i}{b_i}{\mathfrak{p}^s} \quad \text{for all} \ i = 1, \ldots 
, m \ \text{and any} \ s.
\end{equation}
Then \eqref{E:20} implies that $\{ X_s \}$ is a Cauchy sequence in $\mathrm 
O_n(f)_{\mathscr{O}} \subset M_n(\mathscr{K}),$ which therefore converges to 
some $X \in \mathrm O_n(f)_{\mathscr{O}}.$ As in the proof of 
Corollary~\ref{cor:hensels}, we conclude that $\congruence{X}{X_s}{\mathfrak 
p^s},$ so
\[
\congruence{X a_i \equiv X_s a_i}{b_i}{\mathfrak{p}^s} \quad \text{for all} \ 
s,
\]
implying that $X a_i = b_i$ for $i = 1, \dots, m$, as required.

To construct $X_1$  satisfying \eqref{E:20} and \eqref{E:30} for $s=1,$ we 
observe that by applying the usual Witt's theorem to the vector space 
$\bar{\mathscr{L}}$ over $k$ and the reduction of $f$ modulo $\mathfrak{p}$ 
(which is nondegenerate), we can find $Z_1 \in M_n(\mathscr{O})$ such that
\[
\congruence{^tZ_1 F Z_1}{F}{\mathfrak{p}} \quad \text{and} \quad 
\congruence{Z_1 a_i}{b_i}{\mathfrak{p}} \quad \text{for} \ i = 1, \dots, m.
\]
Then by Corollary~\ref{cor:hensels}, there exists $X_1 \in \mathrm 
O_n(f)_{\mathscr{O}}$ satisfying $\congruence{X_1}{Z_1}{\mathfrak{p}}$ and 
possessing thereby the required properties.

Suppose that the matrices $X_1, \ldots , X_s$ have already been constructed. 
Then applying Lemma~\ref{lm:main} to the systems of vectors $\{ X_s a_1, 
\ldots , X_s a_m \}$ and $\{ b_1, \ldots , b_m \}$ in $\mathscr{L}$ yields $X 
\in M_n(\mathscr{O})$ such that
\begin{gather*}
\congruence{X}{E_n}{\mathfrak{p}^s},\\
\congruence{^tXFX}{F}{\mathfrak{p}^{s+1}},
\end{gather*}
and
\[
\congruence{X X_s a_i}{b_i}{\mathfrak{p}^{s+1}}.
\]
Again, by Corollary~\ref{cor:hensels} there exists $X_{s+1} \in \mathrm 
O_n(f)_{\mathscr{O}}$ with the property
\[
\congruence{X_{s+1}}{XX_s}{\mathfrak{p}^{s+1}},
\]
and by our construction such $X_{s+1}$ does satisfy \eqref{E:20} and 
\eqref{E:30} for $s+1,$ completing the proof of the first assertion in $(*)$.

For the second assertion, it is enough to show that if $2m + 1 \leqslant n,$ 
then  $\mathrm O_n(f)_{\mathscr{O}}$ contains a matrix $X$ having determinant 
$-1$ and satisfying $X a_i = a_i$ for all $i = 1, \ldots , m.$ Since the 
reduction $\bar{\mathscr{L}}$ is nondegenerate, there exist $a_{m+1}, \ldots , 
a_{m+r} \in \mathscr{L},$ where $r \leqslant m,$ such that $\bar{a}_1, \ldots, 
\bar{a}_{m+r}$ span a nondegenerate subspace of $\bar{\mathscr{L}}.$ Then the 
lattice $\mathscr{M} = \mathscr{O}a_1 + \cdots + \mathscr{O}a_{m+r}$ is 
unimodular, and therefore $\mathscr{L} = \mathscr{M} \perp 
\mathscr{M}^{\perp}$, where $\mathscr{M}^{\perp}$ is the orthogonal complement 
of $\mathscr{M}$ in $\mathscr{L}.$ Clearly, the lattice $\mathscr{M}^{\perp}$ 
is unimodular, hence contains a vector $c$ such that $f(c) \not\equiv 0 \ 
({\rm mod} \ {\mathfrak p}).$ Then for $X$ one can take (the matrix of) the 
reflection $\tau_c.$
\end{proof}

One corollary of $(*)$ is worth mentioning. Given a nonzero vector $a \in 
\mathscr{L},$ we define its level $\lambda(a)$ (relative to $\mathscr{L}$) as 
follows:
\[
\lambda (a) = \max \{ l \geqslant 0 \mid a \in \mathfrak{p}^l \mathscr{L} \}.
\]

\begin{corollary}\label{cor:witt}
Given two nonzero vectors $a, b \in \mathscr{L}$ such that $f(a)=f(b)$, a 
transformation $X \in \mathrm{O}_n(f)_{\mathscr{O}}$ with the property $X a = 
b$ exists if and only if $\lambda (a) = \lambda (b)$.
\end{corollary}
\begin{proof}
One implication immediately follows from the fact that a matrix in 
$\mathrm{GL}_n(\mathscr{O})$ preserves the level of any vector. For the other 
implication, we write $a$ and $b$ in the form $a= \pi^{\lambda} a_0$, $b = 
\pi^{\lambda} b_0$, where $\lambda = \lambda(a) = \lambda(b).$ Then $f(a_0) = 
f(b_0)$ and the reductions $\bar{a}_0$, $\bar{b}_0 \in \bar{\mathscr{L}}$ are 
nonzero. By $(*)$, there exists $X \in \mathrm O_n(f)_{\mathscr{O}}$ such that 
$Xa_0 = b_0,$ and then also $Xa = b.$
\end{proof}

\begin{remarks}
1. It is worth observing that with some extra work, $(*)$ can be extended to 
discretely valued but not necessarily complete fields $\mathscr{K}.$ Indeed, 
let $\mathscr{K}_v$ be the completion of $\mathscr{K},$ and $\mathscr{O}_v$ be 
the valuation ring in $\mathscr{K}_v.$ Consider the algebraic 
$\mathscr{K}$-group $\mathscr{G} = \mathbf{O}_n(f).$ Now, given $a_1, \ldots , 
a_m$ and $b_1, \ldots , b_m \in \mathscr{L}$ as in $(*)$, we let $\mathscr{H}$ 
denote the stabilizer of all the $a_i$'s in $\mathscr{G}.$ Applying 
respectively the usual Witt's theorem (over $\mathscr{K}$) and $(*)$, we will 
find $X_{\mathscr{K}} \in \mathscr{G}_{\mathscr{K}}$ and $X_v \in 
\mathscr{G}_{\mathscr{O}_v}$ such that
\[
X_{\mathscr{K}} a_i = X_v a_i = b_i \quad \text{for all} \ i = 1, \ldots , m.
\]
Then $Y := X_v^{-1}X_{\mathscr{K}} \in \mathscr{H}_{\mathscr{K}_v}.$ However, 
$\mathscr{H}_{\mathscr{K}}$ is dense in $\mathscr{H}_{\mathscr{K}_v}$ in the 
topology defined by $v.$ (If the subspace spanned by $a_1, \ldots , a_m$ is 
nondegenerate, this is basically proved in 
\cite[Prop.~7.4]{platonov-rapinchuk-book}; the general case is reduced to this 
one by splitting off the unipotent radical of $\mathscr{H}.$) It follows that 
$\mathscr{H}_{\mathscr{K}_v} = \mathscr{H}_{\mathscr{O}_v} 
\mathscr{H}_{\mathscr{K}}.$ Writing $Y$ in the form $Y = Z_v 
Z_{\mathscr{K}}^{-1}$ with $Z_v \in \mathscr{H}_{\mathscr{O}_v}$ and 
$Z_{\mathscr{K}} \in \mathscr{H}_{\mathscr{K}},$ we obtain
\[
X := X_v Z_v = X_{\mathscr{K}} Z_{\mathscr{K}} \in \mathscr{G}_{\mathscr{O}_v} 
\cap \mathscr{G}_{\mathscr{K}} = \mathscr{G}_{\mathscr{O}}
\]
and $X a_i = X_{\mathscr{K}} a_i = b_i$ for all $i = 1, \ldots , m.$

\smallskip

2. As was pointed out by the anonymous referee of the earlier version of the 
paper, the result actually holds for arbitrary local rings and can be derived 
from \cite[Satz~4.3]{kneser} or \cite[Thm.~1.2.2]{kitaoka}.
\end{remarks}

\section{The quadric $Q_s$}\label{sec:quadric}

In this section, we return to the notations and conventions introduced in 
\S\S2--3. To complete the proof of the Main Theorem in the next section, we 
need to figure out when for a given $g \in G_{\mathcal{O}(S)}$ one can choose 
$s , t \in L:= \mathcal{O}(S)e_1 + \cdots + \mathcal{O}(S)e_n$ such that the 
triple $\zeta = (g , s , t)$ belongs to $Z$ and satisfies the assumptions of 
Theorem~\ref{prop:application-witt}. We notice that if $s$ has already been 
chosen so that $(g , s) \in Y$ then the $t$'s for which $(g , s , t)$ belongs 
to $Z$ lie on the following quadric
\[
Q_s = \{ x \in \langle s, a \rangle^{\perp} \mid f(x) = f(b) \}.
\]
So, in this section we will examine some arithmetic properties of $Q_s$ for an 
arbitrary $s \in W$ such that the space $\langle s , a \rangle$ is 
2-dimensional and nondegenerate.
\begin{lemma}\label{lem:aboutqs}
\
\begin{myenumerate}
\item For every $v \in V^K_{\infty}$, $(Q_s)_{K_v} \neq \emptyset$.

\item If $n \geqslant 6$, then $(Q_s)_K \neq \emptyset$.

\item Suppose that $s \in L$. Then for every $v \notin S \cup V_0$, 
$(Q_s)_{\mathcal O_v} \neq \emptyset$.
\end{myenumerate}
\end{lemma}
\begin{proof}
(i) This is obvious if $v$ is complex, so suppose that $v$ is real. Then by 
our construction $n_v^+ \geqslant 3$, implying that the restriction of $f$ to 
$\langle s, a \rangle^{\perp}$ has at least one positive square. Since $f(b) = 
\alpha_{n-1} > 0$ in $K_v = \mathbb{R}$, our assertion follows.

\vskip1mm

(ii) If $n \geqslant 6$ then $\dim \langle s , a \rangle^{\perp} \geqslant 4.$ 
As a nondegenerate  quaternary quadratic form over a (non-Archimedean) local 
field represents every nonzero element \cite[63:18]{omeara}, we conclude that 
$(Q_s)_{K_v} \neq \emptyset$ for all $v \in V^K_f$. Combining this with (i) 
and applying the Hasse--Minkowski theorem \cite[66:4]{omeara}, we obtain our 
claim.

\vskip1mm

(iii) We will show that there is a unimodular ${\mathcal O}_v$-sublattice $M 
\subset L_v := L \otimes_{{\mathcal O}(S)} {\mathcal O}_v$ of rank $\leqslant 
3$ containing $s$ and $a$. Then $L_v = M \perp M^{\perp}$ with $M^{\perp}$ 
unimodular of rank $\geqslant 2.$ Since $f(b) \in {\mathcal O}_v^{\times}$, 
there exists $x \in M^{\perp}$ such that $f(x) = f(b)$ \cite[92:1b]{omeara}, 
and  then $x \in (Q_s)_{\mathcal O_v}$, as required. To construct such an $M$, 
we let  $(s_1, \ldots , s_n)$ denote the coordinates of $s$ in the basis $e_1, 
\ldots , e_n.$ Set $u = s_1e_1 + \cdots  + s_{n-1}e_{n-1}.$ As $u \in L$ and 
$u \neq 0,$ we can write $u = {\pi}_v^d u_0$ where $\pi_v \in {\mathcal O}_v$ 
is a uniformizing element, $d \geqslant 0$ and $u_0 \in L_v \setminus \pi_v 
L_v$. If $f(u_0) \in \mathcal O_v^{\times}$, then in view of $u_0 \perp a,$ 
the sublattice $M = \mathcal O_v a + \mathcal O_v u_0$ is as desired. Now, 
suppose that $f(u_0) \in {\mathfrak p}_v = \pi_v \mathcal O_v$. Since the 
sublattice $N = \mathcal O_v a$ is unimodular we have $L_v = N \perp 
N^{\perp}$ with $N^{\perp}$ unimodular; notice that $u_0 \in N^{\perp}$. The 
reduction $\overline{\left(N^{\perp}\right)}^{(v)} = N^{\perp} 
\otimes_{\mathcal{O}_v} k_v$ being a nondegenerate quadratic space over $k_v = 
\mathcal O_v / \mathfrak p_v,$ one can find $u_1 \in N^{\perp}$ so that the 
images of $u_0$ and $u_1$ in $\overline{\left(N^{\perp}\right)}^{(v)}$ form a 
hyperbolic pair. Then the $\mathcal O_v$-sublattice $M$ with the basis $a,$ 
$u_0$ and $u_1$ is as required.
\end{proof}

\begin{lemma}\label{lem:strongforqs}
\
\begin{myenumerate}
\item Suppose that $n \geqslant 6$. Given  $v \in S,$ there exists an open set 
$\mathcal U_v \subset \mathbf W_{K_v}$ such that $\mathcal{U}_v \cap X \neq 
\emptyset$ and for any $s \in W \cap \mathcal U_v$, the quadric $Q_s$ has 
strong approximation with respect to $S$. Moreover, if there exists $v \in 
V^K_{\infty}$ with the property $i_v(f) \geqslant 2$, then $Q_s$ has strong 
approximation with respect to $S$ for \emph{any} $s$.

\item Suppose that $n = 5$ and $v \in S$ is non-Archimedean. There exists an 
open subset $\mathcal U_v \subset \mathbf W_{K_v}$ with the property 
$\mathcal{U}_v \cap X \neq \emptyset$  such that for $s \in W \cap \mathcal 
U_{v}$ one has $(Q_s)_{K_v} \neq \emptyset$ and moreover if $(Q_s)_K \neq 
\emptyset$ then $Q_s$ has strong approximation with respect to $S.$
\end{myenumerate}
\end{lemma}
\begin{proof}
(i) It follows from the theorem in the Appendix and 
Lemma~\ref{lem:aboutqs}(ii) that  a~necessary and sufficient condition for 
strong approximation in $Q_s$ is that $(Q_s)_S = \prod_{v \in S} (Q_s)_{K_v}$ 
be noncompact. If $v \in V^K_{\infty}$ is such that $i_v(f) \geqslant 2$, then 
for any $s$ the space $\langle s, a \rangle ^{\perp}$ is $K_v$-isotropic and 
therefore $(Q_s)_{K_v}$ is noncompact, hence our second assertion. For the 
first assertion, we observe that the space $\langle a, b \rangle^{\perp}$ is 
$K$-isotropic by our construction, and besides there exists $s_0 \in \langle a 
, b \rangle$ such that $f(s_0) = f(a)$ and $\langle a, s_0 \rangle = \langle 
a, b \rangle.$  The fact that the subgroup of squares $K_v^{\times 2}$ is open 
in $K_v^{\times}$ implies that there exists an open set $\mathcal{U}_v \subset 
\mathbf{W}_{K_v}$ containing $s_0$ such that for any $s \in \mathcal{U}_v,$ 
the spaces $\langle s, a \rangle$ and $\langle s_0, a \rangle$ are isometric 
over $K_v.$  Then it follows from Witt's theorem that the space $\langle s, a 
\rangle^{\perp}$ is $K_v$-isotropic, so the set $\mathcal U_v$ is as required.

\vskip1mm

(ii) Pick $c \in  W$ orthogonal to $a$ and  $b$  so that $f(c) = -f(b)$, and 
let $U$ be the orthogonal complement in $W$ to $a,$ $b,$ and  $c.$  Since 
$\dim U = 2$ and $v$ is non-Archimedean, the set of nonzero values of $f$ on 
$U \otimes_K K_v$ consists of more than one coset modulo $K_v^{\times 2},$ so 
there exists an anisotropic $u \in U$ such that $f(u) \notin -f(c) K_v^{\times 
2},$ and then the space $\langle c , u \rangle$ is $K_v$-anisotropic. Pick $u' 
\in U$ orthogonal to $u$. Then the space $\langle a, u' \rangle^{\perp} = 
\langle b, c, u \rangle$ is $K_v$-isotropic (viz., $f(b + c) = 0$), but the 
space $\langle a, b, u' \rangle^{\perp} = \langle c, u \rangle$ is 
$K_v$-anisotropic. As in the proof of (i), we pick $s_0 \in \langle a , u' 
\rangle$ so that $f(s_0) = f(a)$ and $\langle a, s_0 \rangle = \langle a , u' 
\rangle,$ and then find an open subset $\mathcal{U}_v \subset 
\mathbf{W}_{K_v}$ containing $s_0$ such that for any $s \in \mathcal{U}_v$ the 
subspaces $\langle a , s \rangle$ and $\langle a , s_0 \rangle$ are isometric 
over $K_v.$ If now $s \in W \cap \mathcal{U}_v,$ then it follows from Witt's 
theorem that the space $\langle a , s \rangle^{\perp}$ is $K_v$-isotropic, 
implying not only that $(Q_s)_{K_v} \neq \emptyset$ but in fact also that 
$(Q_s)_{K_v}$ is noncompact. Furthermore, if $d \in (Q_s)_{K_v}$ then the 
space $\langle a , s , d \rangle^{\perp}$ is $K_v$-isometric to $\langle a , 
s_0 , b \rangle^{\perp} = \langle a , b , u' \rangle^{\perp},$ hence 
$K_v$-anisotropic. Thus, if $(Q_s)_K \neq \emptyset,$ then by the theorem in 
the Appendix, $Q_s$ has strong approximation with respect to $S$.
\end{proof}

\section{Proof of the Main Theorem}\label{sec:proof}

For convenience of reference we will list some elementary results about groups 
with bounded generation.

\begin{lemma}\label{prop:abstractbg}
Let $\Gamma$ be a group,  and $\Delta$ be its subgroup.
\begin{myenumerate}
\item If $[\Gamma : \Delta] < \infty$ then bounded generation of $\Gamma$ is 
equivalent to bounded generation of $\Delta.$

\item If $\Gamma$ has {\rm (BG)} then so does any homomorphic image of 
$\Gamma$.

\item If $\Delta \vartriangleleft \Gamma$ and both $\Delta$ and 
$\Gamma/\Delta$ have {\rm (BG)} then $\Gamma$ also has {\rm (BG)}.
\end{myenumerate}
\end{lemma}
\begin{proof}
All these assertions, except for the fact that in (i), (BG) of $\Gamma$ 
implies (BG) of $\Delta,$ immediately follow from the definition. A detailed 
proof of the remaining implication is given, for example, in \cite{murty-95}.
\end{proof}

It follows from Lemma~\ref{prop:abstractbg}(i) that given two commensurable 
subgroups $\Delta_1$ and $\Delta_2$ of $\Gamma$ (which means that $\Delta_1 
\cap \Delta$ has finite index in both $\Delta_1$ and $\Delta_2$), (BG) of one 
of them is equivalent to (BG) of the other. In particular, if $G$ is an 
algebraic group over a number field $K$, then (BG) of one $S$-arithmetic 
subgroup of $G$ implies (BG) of all $S$-arithmetic subgroups of $G.$ 
Furthermore, if $\pi \colon G_1 \to G_2$ is a $K$-defined isogeny of algebraic 
$K$-groups  and $\Gamma$ is an $S$-arithmetic subgroup of $G_1$, then (BG) of 
$\Gamma$ is equivalent to (BG) of $\pi(\Gamma).$ Since the latter is an 
$S$-arithmetic subgroup of $G_2$ (see, for example, 
\cite[Thm.~5.9]{platonov-rapinchuk-book}), we obtain the following.

\begin{lemma}\label{prop:isogeny}
Let $\pi \colon G_1 \to G_2$ be a $K$-defined isogeny of algebraic $K$-groups, 
where $K$ is a number field. Then {\rm (BG)} of one $S$-arithmetic subgroup in 
$G_1$ or $G_2$ implies {\rm (BG)} of all $S$-arithmetic subgroups in $G_1$ and 
$G_2.$
\end{lemma}

Applying this lemma to the universal cover $\mathbf{Spin}_n(f) 
\stackrel{\pi}{\longrightarrow} \mathbf{SO}_n(f),$ we see that to prove the 
Main Theorem it is enough to show that for $G = \mathbf{Spin}_n(f),$ the group 
$G_{\mathcal{O}(S)},$ defined in terms of our fixed realization, is boundedly 
generated. Our argument will use the following simple observation.

\begin{lemma}\label{lm:dobavka}
Let $\Gamma$ be a group, and $\Delta$ be its subgroup of finite index. If 
there exist $\gamma, \gamma_1, \ldots , \gamma_s \in \Gamma$ such that 
$\gamma\Delta \subset \langle \gamma_1 \rangle \cdots \langle \gamma_s 
\rangle$, then $\Gamma$ has {\upshape (BG)}.
\end{lemma}
\begin{proof} Let $x_1, \dots , x_n$ be a system of left coset
representatives for $\Delta$ in  $\Gamma.$  Then
\[
\Gamma \ = \ \bigcup_{i=1}^n x_i \Delta \ \subset \ \bigcup_{i=1}^n x_i 
\langle \gamma \rangle \langle \gamma_1 \rangle \cdots \langle \gamma_s 
\rangle,
\]
implying that  $\Gamma = \langle x_1 \rangle \cdots \langle x_n \rangle 
\langle \gamma \rangle \langle \gamma_1 \rangle \cdots \langle \gamma_s 
\rangle$.
\end{proof}

To proceed with the proof of the Main Theorem, we need to introduce some 
additional notations. For $g \in G,$ the fiber over $g$ of the projection $Y 
\to G$ can (and will) be identified with
\begin{equation}\label{E:100}
B_g = \{ s \in \mathbf W \mid  (s | g(b))= 0, \ f(s) = f(a) \};
\end{equation}
notice that $B_g$ is a quadric in an $(n-1)$-dimensional vector space ($= 
g(b)^{\perp}$). Furthermore, for $v \in V^K$ we denote
\[
\mathcal P_v =
\begin{cases}
(P_0)_{K_v} & \text{if } v \in S,\\
(P_0)_{K_v} \cap  P_{\mathcal O_v} & \text{if } v \notin S,
\end{cases}
\]
where $P_0 \subset  P$ is the Zariski-open set introduced in Lemma~\ref{L:1A}, 
and set $\mathcal G_v = {\mu}(\mathcal P_v)$. It follows from the surjectivity 
of $d_h\mu$ at all points  $h \in P_0$ (see Lemma~\ref{L:1A}) and the Implicit 
Function Theorem \cite[pp.~83--85]{serre} that $\mathcal G_v$ is open in 
$\tg_{K_v}$.

\begin{proposition}\label{prop:witt2}
If $n \geqslant 6$ and $i_K(f) \geqslant 2$ then
\[
\tg_{\oo(S)} \cap \prod_{v \in V_0} \mathcal G_v \subset {\mu}({P}_{\oo(S)}).
\]
\end{proposition}
\begin{proof}
Fix $g \in \tg_{\oo(S)} \cap \prod_{v \in V_0} \mathcal G_v$. Then for each $v 
\in V_0$, one can we pick  $h_v \in \mathcal P_v$ so that ${\mu}(h_v) = g,$ 
hence ${\phi}(h_v) = (g, s_v, t_v)$. It again follows from Lemma \ref{L:1A} 
and the Implicit Function Theorem that the map $(\varepsilon \circ {\phi})_v$ 
is open at $h_v$, implying that one can pick an open  neighborhood $\Sigma_v 
\subset (B_g)_{\oo_v}$ of $s_v$ satisfying
\begin{equation}\label{E:9}
(g, \Sigma_v) \subset \varepsilon({\phi}(\mathcal P_v)).
\end{equation}
Clearly, $g(a) \in B_g,$ in particular, $(B_g)_K \neq \emptyset.$ Furthermore, 
the orthogonal complement of $g(b)$ is isometric to the orthogonal complement 
of $b,$ hence $K$-isotropic, so it follows from \eqref{E:100} that $(B_g)_S$ 
is noncompact. Since $n-1 \geqslant 5,$ by the theorem in the Appendix, $B_g$ 
has strong approximation with respect to $S,$ and therefore one can find
\begin{equation}\label{E:10}
s \in (B_g)_{\oo(S)} \cap \prod_{v \in V_0} \Sigma_v.
\end{equation}
According to Lemma~\ref{lem:strongforqs}(i), the corresponding quadric $Q_s$ 
(see \S \ref{sec:quadric}) has strong approximation with respect to $S$. 
Taking into account that $Q_s = \{ t \mid (g , s , t) \in \varepsilon^{-1}(g , 
s) \}$ and that according to \eqref{E:9} and \eqref{E:10} one has 
$\varepsilon^{-1}(g, s) \cap {\phi}(\mathcal P_v) \neq \emptyset$ for all $v 
\in V_0,$ we conclude that there exists $t$ such that
\begin{equation}\label{E:11}
\zeta := (g, s, t) \in {Z}_{\oo(S)} \cap \prod_{v \in V_0} {\phi}({P}_{\oo_v}).
\end{equation}
Then it follows from Theorem~\ref{prop:application-witt} that 
${\phi}^{-1}(\zeta)_{\oo(S)} \neq \emptyset$, and therefore $g \in {\mu}( 
P_{\oo(S)})$, proving the proposition.
\end{proof}

An analog of Proposition~\ref{prop:witt2} for the case where $i_K(f) = 1$ 
requires a bit more work, especially if $n = 5$.

\begin{proposition}\label{prop:witt1}
Suppose that $n \geqslant 5$, $i_K(f) = 1$, and $S$ contains a non-Archimedean 
valuation. Then
\begin{equation}\label{E:14}
\tg_{\oo(S)} \cap \prod_{v \in V_0} \mathcal G_v  \subset {\mu}({P}_{\oo(S)}) 
{\mu}({P}_{\oo(S)}) {\mu}({P}_{\oo(S)})^{-1}.
\end{equation}
\end{proposition}
\begin{proof}
By our assumption, one can  pick in $S$ an Archimedean valuation $v_1$ and a 
non-Archimedean valuation $v_2.$ Let ${\mathcal U}_{v_2} \subset {\mathbf 
W}_{K_{v_2}}$ be an open subset with the properties described in 
Lemma~\ref{lem:strongforqs}, i.e., ${\mathcal U}_{v_2} \cap X \neq \emptyset$ 
and for any $s \in {\mathcal U}_{v_2} \cap X_K,$ the quadric $Q_s$ has strong 
approximation with respect to $S$ if either $n \geqslant 6$ or $n = 5$ and 
$(Q_s)_K \neq \emptyset;$ in addition, for such $s$ one can guarantee that 
$(Q_s)_{K_v} \neq \emptyset$ if $n = 5.$ It now follows from Lemma~\ref{L:2A} 
that for the map $\eta$ introduced therein, the set
\[
\mathcal P'_{v_2} = \eta^{-1}(\mathcal U_{v_2} \cap X_{K_{v_2}}) \cap ( 
P_0)_{K_{v_2}}
\]
is a nonempty open subset of ${P}_{K_{v_2}}$. Then  as above we conclude that 
$\mathcal G'_{v_2} := {\mu}(\mathcal P'_{v_2})$ is a nonempty open subset of 
$\tg_{K_{v_2}}$.

By strong approximation, ${P}_{\oo(S)}$ is dense in ${P}_{S \setminus \{ v_1 
\}} \times \prod_{v \in V_0} {P}_{\oo_v}$, which in view of 
Proposition~\ref{prop:phipsi} implies that the closure of 
${\mu}({P}_{\oo(S)})$  in ${G}_{(S \cup V_0) \setminus \{ v_1 \}}$ contains 
$G_{S \setminus \{ v_1 \}} \times \prod_{v \in V_0} {\mathcal G}_v.$ Since the 
${\mathcal G}_v$'s are open, we conclude that the closure of ${\mathcal B} := 
{\mu}({P}_{\oo(S)}) {\mu}({P}_{\oo(S)})^{-1}$ in ${G}_{(S \cup V_0) \setminus 
\{ v_1 \}}$ contains $G_{S \setminus \{ v_1 \}} \times \prod_{v \in V_0} 
{\mathcal E}_v$ for some open neighborhoods of the identity ${\mathcal E}_v 
\subset G_{K_v},$ $v \in V_0.$ It follows that given an element $g$ belonging 
to the left-hand side of the inclusion (\ref{E:14}), there exists $h \in 
{\mathcal B}^{-1}$ such that
\[
gh \in \left( \prod_{v \in (V_0 \cup S) \setminus \{ v_1 , v_2 \}} \mathcal G_v 
\right) \times  \mathcal G'_{v_2}.
\]
Thus, it is enough to show that
\begin{equation}\label{E:12}
\tg_{\oo(S)} \bigcap \left[ \left( \prod_{v \in (V_0 \cup S) \setminus \{ v_1 , 
v_2 \}} \mathcal G_v \right) \times  \mathcal G'_{v_2} \right] \subset 
{\mu}({P}_{\oo(S)}).
\end{equation}
Fix a $g$ belonging to the left-hand side of the inclusion \eqref{E:12}. As in 
the proof of Proposition~\ref{prop:witt2}, we can find open sets  $\Sigma_v 
\subset (B_g)_{\oo_v}$  for $v \in (S \cup V_0) \setminus \{ v_1 , v_2 \}$ 
such that
\[
(g, \Sigma_v) \subset \varepsilon({\phi}(\mathcal P_v))
\]
and also an open set $\Sigma'_{v_2} \subset (B_g)_{K_{v_2}}$ such that
\[
(g, \Sigma'_{v_2}) \subset \varepsilon({\phi}(\mathcal P'_{v_2})).
\]
As in the proof of Proposition~\ref{prop:witt2}, we use the theorem in the 
Appendix to conclude that $B_g$ has strong approximation with respect to $\{ 
v_1 \}$, so one can find
\[
s \in (B_g)_{\oo(S)} \bigcap \left[ \left( \prod_{v \in (S \cup V_0) \setminus 
\{ v_1 , v_2\}} \Sigma_v \right) \times \Sigma'_{v_2} \right].
\]
Then for each $v \in (S \cup V_0) \setminus \{ v_1 \}$, we have 
$\varepsilon^{-1}(g , s)_{K_v} \neq \emptyset$ implying that $(Q_s)_{K_v} \neq 
\emptyset.$ Furthermore, the non-emptiness of $(Q_s)_{K_v}$ for $v = v_1$ 
follows from Lemma~\ref{lem:aboutqs}(i) as $v_1$ is Archimedean, and for $v 
\notin S \cup V_0$ --- from Lemma~\ref{lem:aboutqs}(iii) as $s \in L.$  So, by 
the Hasse--Minkowski theorem \cite[66:4]{omeara}, $(Q_s)_K \neq \emptyset$. 
Since $s \in {\mathcal U}_{v_2},$ by Lemma~\ref{lem:strongforqs}, $Q_s$ has 
strong approximation with respect to $S$ in all cases. The rest of the 
argument repeats verbatim the corresponding part of the proof of 
Proposition~\ref{prop:witt2}: we use strong approximation for $Q_s$ to find a 
$t$ for which the triple $\zeta =(g, s, t)$ satisfies \eqref{E:11}; then by 
Theorem~\ref{prop:application-witt}, ${\phi}^{-1}(\zeta)_{\oo(S)} \neq 
\emptyset$. This implies that $g \in {\mu}({P}_{\oo(S)})$, and the proposition 
follows.
\end{proof}

\begin{proof}[Proof of the Main Theorem]
As we explained in the beginning of this section, it is enough to establish 
bounded generation of $G_{\oo(S)}.$ For this, we will argue by induction on 
$n$. First, we will consider the case $i_K(f) \geqslant 2$. In this case we 
can assume without any loss of generality that the basis $e_1, \ldots , e_n$ 
is chosen so that the space spanned by $e_1, \ldots, e_4$ has Witt index two. 
If $n = 5$,  then the group ${G}$ is $K$-split, so bounded generation of 
${G}_{\oo(S)}$ is a result of Tavgen \cite{tavgen-90}. For $n \geqslant 6$, it 
follows from Proposition~\ref{prop:witt2} that $\mu(P_{\oo(S)})$ contains an 
open subset of ${G}_{\oo(S)}$, and since congruence subgroups form a base of 
neighborhoods of the identity, there exists a congruence subgroup $\Delta 
\subset {G}_{\oo(S)}$ and an element $h \in {G}_{\oo(S)}$ such that
\[
h \Delta \subset \mu(P_{\oo(S)}) = \tg(a)_{\oo(S)} \tg(b)_{\oo(S)} 
\tg(a)_{\oo(S)} \tg(b)_{\oo(S)}.
\]
Since both ${G}(a)_{\oo(S)}$ and ${G}(b)_{\oo(S)}$ are boundedly generated by 
induction hypothesis and $\Delta$ has finite index in ${G}_{\oo(S)}$, bounded 
generation of the latter follows from Lemma~\ref{lm:dobavka}.

Now, suppose that $i_K(f) = 1$ but $S$ contains a non-Archimedean valuation. 
Here the induction starts with $n = 4,$ in which case $G$ is known to be 
$K$-isomorphic to either $\mathbf{SL}_2 \times \mathbf{SL}_2$ or 
$R_{E/K}(\mathbf{SL}_2)$ for a suitable quadratic extension $E/K$ (see, for 
example, \cite[Thms.~5.21 and 5.22]{artin}). In either case, since $S$ contains 
a non-Archimedean valuation, bounded generation of ${G}_{\oo(S)}$ follows from 
bounded generation of $\mathrm{SL}_2(A)$, where $A$ is a ring of $S$-integers 
in a number field having infinitely many units 
\cite{cooke-w-75,carter-keller-preprint}. For $n \geqslant 5$, the argument is 
completed as above using Proposition~\ref{prop:witt1} instead of 
Proposition~\ref{prop:witt2}.
\end{proof}

\begin{appendix}
\section{Appendix}

The purpose of this appendix is to formulate and prove the result on strong 
approximation in quadrics that was used in the proof of the Main Theorem. Let 
$q = q(x_1, \ldots , x_m)$ be a nondegenerate quadratic form in $m \geqslant 3$ 
variables over a number field $K,$ and $Q$ be a quadric given by the equation 
$q(x_1, \ldots , x_m) = a$ where $a \in K^{\times}.$ Fix a nonempty subset $S$ 
of $V^K.$

\begin{named}{Theorem}
Assume that $Q_K \neq \emptyset$ and $Q_S := \prod_{v \in S} Q_{K_v}$ is 
noncompact.
\begin{myenumerate}

\item If $m \geqslant 4$ then $Q$ has strong approximation with respect to $S.$

\item If $m = 3$ then $Q$ has strong approximation with respect to $S$ if and 
only if the following condition holds: Let $x \in Q_K$ and let $g$ be the 
restriction of $q$ to the orthogonal complement of $x$ in $K^3;$ then either 
$g$ is $K$-isotropic, or $g$ is $K$-anisotropic and there exists $v \in S$  for 
which $g$ is $K_v$-anisotropic and additionally $q$ is $K_v$-isotropic if $v$ 
is real.
\end{myenumerate}
\end{named}

Assertion (i) is proved, for example, in \cite[104:3]{omeara}, where it is then 
used to establish strong approximation for $\mathbf{Spin}_m(q).$ We have not 
found, however, a proof of assertion (ii) in the literature. As was pointed out 
in \cite{rapin-88}, both facts can be derived from the analysis of strong 
approximation in the homogeneous spaces $G/H$ which relies on the strong 
approximation theorem for algebraic groups and results on Galois cohomology. 
Such analysis for the cases where $G$ is a connected simply connected $K$-group 
and $H$ is either its connected simply connected $K$-subgroup or a $K$-subtorus 
(which are sufficient for the proof of the theorem) was given in 
\cite{rapin-88}; the case of an arbitrary reductive $H$ was independently 
considered in \cite{borovoi-89}. In our exposition we will follow 
\cite{rapin-90}.

First, we establish the following criterion of strong approximation which 
easily translates into the  language of Galois cohomology.

\begin{lemma}\label{AL:1}
Let $X = G/H$, where $G$ is a connected $K$-group and $H$ is its connected 
$K$-subgroup. If $G$ has strong approximation with respect to $S$ then the 
closure of $X_K$ in $X_{A(S)}$ coincides with $G_{A(S)} X_K = \{ gx \mid  g \in 
G_{A(S)}, \ x \in X_K \}.$ Thus, $X$ has strong approximation with respect to 
$S$ if and only if the map of the orbit spaces $G_K \backslash X_K \to G_{A(S)} 
\backslash X_{A(S)}$ is surjective.
\end{lemma}
\begin{proof}
It follows from the Implicit Function Theorem that for every $v \in V^K$ and 
any $x_v \in X_{K_v},$ the orbit $G_{K_v} x_v$ is open in $X_{K_v}$ 
\cite[\S3.1, Cor.~2]{platonov-rapinchuk-book}. Moreover, for almost all $v \in 
V^K_f,$ the group $G_{{\mathcal O}_v}$ acts on $X_{{\mathcal O}_v}$ 
transitively (this is a consequence of Hensel's lemma and the fact that for 
almost all $v$ there exist smooth (irreducible) reductions 
$\underline{G}^{(v)}$, $\underline{H}^{(v)}$ and $\underline{X}^{(v)} = 
\underline{G}^{(v)}/\underline{H}^{(v)}$, so $\underline{G}^{(v)}$ acts on 
$\underline{X}^{(v)}$ transitively by Lang's theorem (see, for example, 
\cite[\S16]{Borel}). Thus,  for any $x \in X_{A(S)},$ the orbit $G_{A(S)}x$ is 
open in $X_{A(S)}.$ We conclude that the complement of $G_{A(S)}X_K$ in 
$X_{A(S)}$ is open, hence $G_{A(S)}X_K$ is a closed subset of $X_{A(S)}$ 
containing $X_K.$ On the other hand, strong approximation in $G$ implies that 
$X_K = G_K X_K$ is dense in $G_{A(S)} X_K,$ and all our assertions follow.
\end{proof}

To give a cohomological interpretation, we recall that for any field extension 
$P/K$, there is a natural bijection
\[
G_P \backslash X_P \ \simeq \ \Ker \left( H^1 (P, H) \to H^1 (P, G) \right)
\]
(see \cite{galoiscoh} for the details and unexplained notations). In the adelic 
setting, for any finite Galois extension $L/K,$ there is a bijection
\begin{multline*}
G_{A(S)} \backslash (X_{A(S)} \cap \alpha(G_{A(S) \otimes_K L})) \ \simeq \\
\Ker \left( H^1 (L/K, H_{A(S) \otimes_K L}) \to H^1 (L/K, G_{A(S) \otimes_K L}) 
\right)
\end{multline*}
where $\alpha \colon G \to G/H = X$ is the canonical map. Given an algebraic 
$K$-group $D$, we let $H^1 (K, D)_{A(S)}$ denote the direct limit of the sets 
$H^1 (L/K, D_{A(S) \otimes_K L})$ taken over all finite Galois extensions 
$L/K$; we notice that if $D$ is connected then for a fixed $L/K$ the set $H^1 
(L_w / K_v, D_{{\mathcal O}(L_w)})$ is trivial for almost all $v \in V^K_f$, 
where $w \vert v$, so $H^1 (K, D)_{A(S)}$ can be identified with the set 
$\prod_{v \notin S}' H^1 (K_v, D)$ consisting of $(c_v) \in \prod_{v \notin S} 
H^1 (K_v, D)$ such that $c_v$ is trivial for almost all $v$ (see 
\cite[\S6.2]{platonov-rapinchuk-book}). With these notations, there is a 
bijection
\[
G_{A(S)} \backslash X_{A(S)} \ \simeq \ \Ker \left( H^1 (K, H)_{A(S)} \to H^1 
(K, G)_{A(S)} \right)
\]
Now we can reformulate Lemma~\ref{AL:1}  as follows.

\begin{corollary}\label{AC:1}
Let $X = G/H$ as above. Assume that  $G$ has strong approximation with respect 
to $S.$ Then $X$ has strong approximation with respect to $S$ if and only if 
the natural map
\[
\Ker \left( H^1(K, H) \to H^1 (K, G) \right) \ \longrightarrow \Ker \left( 
H^1(K, H)_{A(S)} \to H^1(K, G)_{A(S)} \right)
\]
is surjective.
\end{corollary}

We recall that to have strong approximation with respect to a finite $S,$ an 
algebraic group $G$ must be connected and simply connected 
\cite[\S7.4]{platonov-rapinchuk-book}, so we will assume that this is the case 
in the rest of this appendix. The cohomological criterion of 
Corollary~\ref{AC:1} immediately leads to the following.

\begin{proposition}\label{AP:1}
Let $X = G/H$ where $G$ has strong approximation with respect to $S$. If $H$ is 
connected and simply connected then $X$ also has strong approximation with 
respect to $S$.
\end{proposition}
\begin{proof}
Since $G$ and $H$ are both simply connected, $H^1 (K_v, G)$ and $H^1 (K_v, H)$ 
are trivial for all $v \in V^K_f$ \cite[Thm~6.4]{platonov-rapinchuk-book}. This 
means that
\begin{multline*}
\Ker \left( H^1 (K, H)_{A(S)} \to H^1 (K, G)_{A(S)} \right) \\
= \prod_{v \in V^K_{\infty} \setminus (V^K_{\infty} \cap S)} \Ker \left( H^1 
(K_v, H) \to H^1 (K_v, G) \right).
\end{multline*}
So, the proposition follows from Corollary~\ref{AC:1} and the fact that the map
\[
\psi \colon \Ker \left( H^1(K, H) \to H^1 (K, G) \right) \longrightarrow 
\prod_{v \in V^K} \Ker \left( H^1 (K_v, H) \to H^1 (K_v, G) \right)
\]
is surjective. This is in fact true for any connected $H.$ Indeed, we have the 
following commutative  diagram:
\[
\begin{CD}
H^1 (K, H) @>>> H^1 (K, G)\\
@V\beta VV @VV\gamma V\\
\prod_{v \in V^K_{\infty}} H^1 (K_v, H) @>>> \prod_{v \in V^K_{\infty}} H^1 
(K_v, G)
\end{CD}
\]
Since $\beta$ is surjective \cite[Prop.~6.17]{platonov-rapinchuk-book} and 
$\gamma$ is injective (``Hasse principle'', 
\cite[Thm.~6.6]{platonov-rapinchuk-book}), the surjectivity of $\psi$ follows.
\end{proof}

Proposition~\ref{AP:1} readily yields assertion (i) of the theorem. Indeed, it 
follows from Witt's theorem that $Q$ is a homogeneous space of $G = 
\mathbf{Spin}_m(q)$ so that  if $x \in Q_K$ then $Q$ can be identified with the 
homogeneous space $X = G/H$ where $H = G(x)$. Clearly, $H = 
\mathbf{Spin}_{m-1}(g)$, where $g$ is the restriction of $q$ to the orthogonal 
complement of $x$; in particular $H$ is connected and simply connected for $m 
\geqslant 4$. Since $Q_S$ is noncompact, $G_S$ is also noncompact, and hence 
has strong approximation with respect to $S$. Thus, strong approximation for $Q 
\cong X$ follows from Proposition~\ref{AP:1}. If $m = 3$ then $H = 
\mathbf{Spin}_2(g)$ is a 1-dimensional torus, so to handle this case we need to 
analyze the cohomological criterion of Corollary~\ref{AC:1} in the situation 
where $H = T$ is a $K$-torus.

So, let $T$ be a $K$-torus of a connected simply connected $K$-group $G$. Fix a 
finite Galois extension $L/K$ that splits $T$. It follows from Hilbert's 
Theorem 90 that
\[
H^1 (K, T) = H^1 (L/K, T)  \quad \text{and} \quad  H^1 (K, T)_{A(S)} = H^1 
(L/K, T_{A(S) \otimes_K L}).
\]
So, the map in Corollary~\ref{AC:1} reduces to the following
\begin{multline*}
\phi \colon \Ker \left( H^1( L/K, T) \to H^1 ( L/K, G) \right)
\longrightarrow\\
\Ker \left( H^1( L/K, T_{A(S) \otimes_K L}) \to H^1 ( L/K, G_{A(S) \otimes_K 
L}) \right).
\end{multline*}
We now let $A$ denote the (full) adelic ring of $K$. It follows from the Hasse 
principle for $G$ that the map $\gamma$ in the following commutative diagram
\[
\begin{CD}
H^1 (L/K, T) @>\alpha>> H^1 (L/K, G)\\
@V\beta VV @VV\gamma V\\
H^1 (L/K, T_{A \otimes_K L}) @>\delta>> H^1 (L/K, G_{A \otimes_K L})
\end{CD}
\]
is injective, so
\begin{equation}\label{E:100}
\beta ( \Ker \alpha) = \Image \beta \cap \Ker \delta.
\end{equation}
Let $C_L(T) = T_{A \otimes_K L}/T_L$ denote the group of classes of adeles of 
$T$ over $L$. The exact sequence
\[
1 \longrightarrow T_L \longrightarrow T_{A \otimes_K L} \longrightarrow C_L(T) 
\longrightarrow 1
\]
gives rise to the exact cohomological sequence
\begin{equation}\label{E:101}
H^1 (L/K, T) \stackrel{\beta}{\longrightarrow} H^1 (L/K, T_{A \otimes_K L}) 
\stackrel{\rho}{\longrightarrow} H^1 (L/K, C_L(T)).
\end{equation}
Writing $A =  A(S) \times K_S$ where $K_S = \prod_{v \in S} K_v$ and using 
\eqref{E:100} in conjunction with the exactness of \eqref{E:101}, we obtain
\begin{multline}\label{E:102}
\Image \phi = \{ x \in \Ker \left( H^1 (L/K, T_{A(S) \otimes_K L}) \to H^1
(L/K, G_{A(S) \otimes_K L}) \right) \mid \text{ there is } \\
y \in \Ker \left( H^1 (L/K, T_{K_S \otimes_K L}) \to H^1 (L/K, G_{K_S \otimes_K 
L}) \right) \text{ with } \rho(x, y) = 0 \}.
\end{multline}
Now we are in a position to give a criterion for strong approximation in $X = 
G/T$ in terms of properties of the map $\rho$.

\begin{proposition}\label{AP:2}
Let $X = G/T$ where $G$ is a simply connected $K$-group and $T$ is a 
$K$-subtorus of $G$. Assume that $G$ has strong approximation with respect to 
$S$. Then $X$ has strong approximation with respect to $S$ if and only if
\begin{multline}\label{E:103}
\rho \left( H^1 (L/K, T_{A(S) \otimes_K L}) \right)\\
\subset \rho \left( \Ker \left( H^1 (L/K, T_{K_S \otimes_K L}) \to H^1 (L/K, 
G_{K_S \otimes_K L}) \right) \right).
\end{multline}
\end{proposition}
\begin{proof}
It follows from \eqref{E:102} and Corollary~\ref{AC:1} that all we need to 
prove is the equality
\begin{multline}\label{E:104}
\rho \left( \Ker \left( H^1 (L/K, T_{A(S) \otimes_K L}) \to H^1 (L/K, G_{A(S)
\otimes_K L}) \right) \right) \\
= \rho \left( H^1 (L/K, T_{A(S) \otimes_K L}) \right).
\end{multline}
Notice that for any  $v \in V^K_f$ and its extension $w \in V^L_f,$ the first 
cohomology $H^1 (L/K, G_{K_v \otimes_K L}) = H^1 (L_w/K_v, G_{L_w})$ is 
trivial. This implies that
\begin{multline*}
\Ker \left( H^1 (L/K, T_{A(S) \otimes_K L}) \to H^1 (L/K, G_{A(S) \otimes_K L})
\right) =\\
\Ker \left( H^1 (L/K, T_{K_{S_{\infty}} \otimes_K L}) \to H^1 (L/K, 
G_{K_{S_{\infty}} \otimes_K L}) \right) \times H^1 (L/K, T_{A(S \cup 
S_{\infty}) \otimes_K L})
\end{multline*}
where $S_{\infty} = V^K_{\infty} \setminus (V^K_{\infty} \cap S)$. Thus, to 
establish \eqref{E:104} it suffices to show that for any $v_0 \in V^K_{\infty}$ 
there exists $v \notin S \cup V^K_{\infty}$ such that
\begin{equation}\label{E:105}
\rho (H^1 (L/K, T_{K_{v_0} \otimes_K L})) = \rho (H^1 (L/K, T_{K_v \otimes_K 
L})).
\end{equation}
Let $X_*(T)$ be the group of cocharacters of $T$ (i.e., $X_*(T) = 
\mathrm{Hom}\; (\mathbf G_m, T))$. It follows from the Nakayama--Tate Theorem 
\cite{voskres} that one can identify
\[
H^1 (L/K, C_L(T)) \quad \text{with} \quad \hat{H}^{-1} (L/K, X_*(T))
\]
and
\[
H^1 (L/K, T_{K_{v_0} \otimes_K L}) = H^1 (L_{w_0}/K_{v_0}, T_{L_{w_0}}) \quad 
\text{with} \quad \hat{H}^{-1} (L_{w_0}/K_{v_0}, X_*(T))
\]
and under these identifications the left-hand side of \eqref{E:105} coincides 
with the image of the corestriction map
\[
\mathrm{Cor}^{\mathrm{Gal}(L/K)}_{\mathrm{Gal}(L_{w_0}/K_{v_0})} \colon 
\hat{H}^{-1} (L_{w_0}/K_{v_0}, X_*(T)) \longrightarrow \hat{H}^{-1} (L/K, 
X_*(T)).
\]
Similarly, the right-hand side of \eqref{E:105} coincides with the image of
\[
\mathrm{Cor}^{\mathrm{Gal}(L/K)}_{\mathrm{Gal}(L_{w}/K_{v})} \colon 
\hat{H}^{-1} (L_{w}/K_{v}, X_*(T)) \longrightarrow \hat{H}^{-1} (L/K, X_*(T))
\]
(we fix extensions $w_0 \vert v_0$ and $w \vert v$). Thus, \eqref{E:105} 
definitely holds if $\mathrm{Gal}(L_{w_0}/K_{v_0}) = 
\mathrm{Gal}(L_{w}/K_{v})$. But for $v_0 \in V^K_{\infty}$, the Galois group 
$\mathrm{Gal}(L_{w_0}/K_{v_0})$ is cyclic, so the existence of $v \notin S \cup 
V^K_{\infty}$ with the same Galois group $\mathrm{Gal}(L_w/K_v)$ follows from 
Chebotarev Density Theorem (see, for example, \cite[Ch.~VII, 
Thm.~13.4]{neukirch}).
\end{proof}

We can now complete the proof of assertion (ii) of the theorem. As we pointed 
out earlier, here $Q$ can be identified with the homogeneous space $X = G/T$, 
where $G = \mathbf{Spin}_3 (q)$ and $T$ is the 1-dimensional torus 
$\mathbf{Spin}_2 (g)$ where $g$ is the restriction of $q$ to the orthogonal 
complement of a chosen point $x \in Q_K$. If $g$ is $K$-isotropic then $T$ 
splits over $L = K$, so \eqref{E:103} trivially holds, and 
Proposition~\ref{AP:2} yields strong approximation in $Q \cong X$.

Suppose now that $g$ is $K$-anisotropic. Then $T$ splits over a quadratic 
extension $L/K,$ with the nontrivial element of $\mathrm{Gal} (L/K)$ acting on 
$X_*(T) \cong \mathbb Z$ as multiplication by $-1$, so
\[
H^1 (L/K, C_L(T)) \cong \hat{H}^{-1} (L/K, X_*(T)) \cong \mathbb Z/ 2 \mathbb 
Z.
\]
Furthermore, by Chebotarev Density Theorem there exists $v \notin S \cup 
V^K_{\infty}$ such that $L_w/K_v$ is a quadratic extension, and then
\[
\rho \left( H^1 (L_w/K_v, T) \right) = H^1 (L/K, C_L(T))
\]
implying that
\[
\rho \left( H^1 (L/K, T_{A(S) \otimes_K L}) \right) = H^1 (L/K, C_L(T)).
\]
Thus, the condition \eqref{E:103} that gives a criterion for strong 
approximation in $X$ boils down to the equality
\[
\rho \left( \Ker \left( H^1 (L/K, T_{K_S \otimes_K L}) \to H^1 (L/K, G_{K_S 
\otimes_K L}) \right) \right) = H^1 (L/K, C_L(T)),
\]
which in turn holds if and only if there is $v \in S$ such that
\begin{equation}\label{E:106}
\Ker \left( H^1 (L_w/K_v, T) \to H^1 (L_w/K_v, G) \right) \neq \{ 1 \}.
\end{equation}
Clearly, \eqref{E:106} holds if $L_w/K_v$ is a quadratic extension (i.e., $g$ 
is $K_v$-anisotropic) and $H^1 (L_w/K_v, G) = \{ 1 \}$ which happens if either 
$v \in V^K_f$ or $q$ is $K_v$-isotropic (notice that in the latter case $G 
\cong \mathbf{SL}_2$ over $K_v$). This proves the presence of strong 
approximation in all cases listed in (ii). It remains to show that in all other 
situations strong approximation does not hold, i.e., \eqref{E:106} fails for 
all $v \in S$. If $T$ splits over $K_v$ then $H^1 (L_w/K_v, T) = \{ 1 \}$, so 
\eqref{E:106} cannot possibly hold. In the remaining case, $v$ is real and $G$ 
is $K_v$-anisotropic. Then $G = \mathbf{SL}_1(\mathbb H)$ where $\mathbb H$ is 
the algebra of Hamiltonian quaternions and $T$ corresponds to a maximal 
subfield of $\mathbb H$. A simple computation shows that the map $H^1 (\mathbb 
C/ \mathbb R, T) \to H^1 (\mathbb C/ \mathbb R, G)$ is a bijection, so again 
\eqref{E:106} fails. (Thus, the 2-dimensional quadric over $\mathbb Q$ given by 
the equation $x_1^2 + x_2^2 - 2x_3^2 = 1$ does not have strong approximation 
with respect to $S = V^{\mathbb Q}_{\infty}$.)

\end{appendix}

%
\newlength{\odin}
\settowidth{\odin}{University of North Carolina at Greensboro}

\newlength{\dva}
\settowidth{\dva}{Department of Mathematics}

\bigskip

\noindent \parbox[t]{\odin}{Department of Mathematical Sciences\\
University of North Carolina at Greensboro\\
Greensboro NC 27402\\
E-mail: igor@uncg.edu} \hfill
\parbox[t]{\dva}{Department of Mathematics\\
University of Virginia\\
Charlottesville VA 22904\\
E-mail: asr3x@virginia.edu}
%

\end{document}